\documentclass[12pt]{amsart}

\usepackage[english]{babel}
\usepackage[utf8x]{inputenc}
\usepackage[T1]{fontenc}
\usepackage{comment}
\usepackage[usenames,dvipsnames]{xcolor} 

\usepackage{mathpazo}
\usepackage[colorlinks, linkcolor=BlueViolet, urlcolor=BlueViolet, citecolor=BlueViolet]{hyperref}
\usepackage{bm}
\usepackage{verbatim}
\usepackage[margin=1in]{geometry}
\usepackage{array,multirow}
\usepackage{graphicx}
\usepackage{tikz-cd}
\usepackage{colonequals}
\usepackage{adjustbox}
\usepackage{booktabs}
\usepackage{comment}

\usepackage[all]{xy}
\usepackage{array}
\usepackage{tikz-cd}

\usepackage{caption}
\usepackage{subcaption}

\usepackage{amssymb,amsmath,latexsym,amsthm}
\usepackage{multirow, cellspace}

\newtheorem{theorem}{Theorem}[section]
\newtheorem{lemma}[theorem]{Lemma}
\newtheorem{proposition}[theorem]{Proposition}

\newtheorem*{theorem*}{Theorem A}
\newtheorem*{theorem'}{Theorem B}
\newtheorem*{theorem"}{Theorem C}

\newtheorem{corollary}[theorem]{Corollary}
\newtheorem{predefinition}[theorem]{Definition}
\newenvironment{definition}{\begin{predefinition}\rm}{\end{predefinition}}
\newtheorem{preremark}[theorem]{Remark}
\newenvironment{remark}{\begin{preremark}\rm}{\end{preremark}}
\newtheorem{prenotation}[theorem]{Notation}

\newtheorem{preexample}[theorem]{Example}
\newenvironment{example}{\begin{preexample}\rm}{\end{preexample}}
\newtheorem{preclaim}[theorem]{Claim}

\newtheorem{prequestion}[theorem]{Question}

\newtheorem{preapplication}[theorem]{Application}

\numberwithin{equation}{section}

\usepackage{amsfonts,amssymb,amsmath,amsthm,latexsym,dsfont}
\usepackage{pifont}
\usepackage{setspace}
\usepackage{mathtools}
\usepackage{epsfig}
\usepackage{amsfonts}
\usepackage{amsmath, amsthm, amssymb}
\usepackage{graphicx}
\usepackage{euscript}
\usepackage{enumerate}
\usepackage{tikz-cd}
\usepackage{listings}
\usepackage{pdfpages}
\usepackage{chngcntr}
\usepackage{pgfplots}
\pgfplotsset{compat=1.17}
\lstdefinelanguage{Magma}%
  {otherkeywords={:=,+:=,-:=,*:=},%
   procnamekeys={function,func,intrinsic,procedure,proc},%
   morekeywords={true,false},%
   morekeywords=[2]{adj,and,cat,cmpeq,cmpne,diff,div,eq,ge,gt,in,is,join,le,lt,%
          meet,mod,ne,notadj,notin,notsubset,or,sdiff,subset,xor},%
   morekeywords=[3]{assigned,break,by,case,catch,continue,declare,default,%
          delete,do,elif,else,end,eval,exists,exit,for,forall,fprintf,if,local,%
          not,print,printf,quit,random,read,readi,repeat,restore,save,select,%
          then,time,to,try,until,vprint,vprintf,vtime,when,where,while},%
   morekeywords=[4]{clear,forward,freeze,iload,import,load},%
   morekeywords=[5]{assert,assert2,assert3,error,require,requirege,requirerange},%
   morekeywords=[6]{car,comp,cop,elt,ext,frac,hom,ideal,iso,lideal,loc,map,%
          ncl,pmap,quo,rec,recformat,rep,rideal,sub},%
   morekeywords=[7]{AbelianGroup,AdditiveCode,AffineAlgebra,Algebra,%
          AssociativeAlgebra,Character,CliffordAlgebra,Design,Digraph,%
          ExtensionField,FPAlgebra,FiniteAffinePlane,FiniteProjectivePlane,%
          Graph,Group,GroupAlgebra,IncidenceStructure,LieAlgebra,LinearCode,%
          LinearSpace,MatrixAlgebra,MatrixGroup,MatrixRing,Monoid,%
          MultiDigraph,MultiGraph,NearLinearSpace,Network,PartialMap,%
          PermutationGroup,PolycyclicGroup,QuaternionAlgebra,Semigroup,%
          ZModule},%
   morekeywords={[8]function,func,intrinsic,procedure,proc,return},%
      sensitive,%
      morecomment=[l]//,%
      morecomment=[s]{/*}{*/},%
      morecomment=[s]{\{}{\}},%
      morestring=[b]"%
  }[keywords,procnames,comments,strings]%
\numberwithin{equation}{section}
\numberwithin{figure}{section}
\numberwithin{table}{section}

\makeatletter
\def\keywords{\xdef\@thefnmark{}\@footnotetext}
\makeatother

\def \PP {{\mathbb P}}

\def \Z {{\mathbb Z}}
\def \cO {{\mathcal O}}
\def \cF {{\mathcal F}}
\def \G {{\mathcal G}}

\def \d {\delta}

\def \t {\tau}

\def \om {\omega}
\def \Om {\Omega}
\def \ph {\varphi}

\def \im {\textup{im}}

\def \ker {\textup{ker }}
\def \im {\textup{Im }}
\def \Gal {\textup{Gal}}
\def \supp {\textup{supp}}
\def \ord {\textup{ord}}
\def \coef {\textup{coef}}
\def \spn {\textup{span}}
\def \tg {\Tilde{g}}
\def \max {\text{max}}
\def \Jac {\text{Jac}}

\def\multiset#1#2{\ensuremath{\left(\kern-.3em\left(\genfrac{}{}{0pt}{}{#1}{#2}\right)\kern-.3em\right)}}

\title{Powers of the Cartier operator on Artin-Schreier covers}
\date{2022}
\author{Steven R. Groen}
\address{Department of Mathematics, University of Warwick, Coventry, CV4 7AL, United Kingdom}
\email{steven.groen@warwick.ac.uk}

\begin{document}
\newcounter{longlist}

\begin{abstract}
For a curve in positive characteristic, the Cartier operator acts on the vector space of its regular differentials. The $a$-number is defined to be the dimension of the kernel of the Cartier operator. In \cite{BoCaASc}, Booher and Cais use a sheaf-theoretic approach to give bounds on the $a$-numbers of Artin-Schreier covers. In this paper, I generalize that approach to arbitrary powers of the Cartier operator, yielding bounds for the dimension of the kernel. These bounds give new restrictions on the Ekedahl-Oort type of Artin-Schreier covers.
\end{abstract}

\maketitle 


\section{Introduction}

Let $k$ be an algebraically closed field of characteristic $p>0$ and let $\pi: Y \to X$ be a finite morphism of smooth, projective, geometrically connected curves over $k$ that is generically Galois with Galois group $G$, branched at a finite set $S \subset X(k)$. In general we wish to express invariants of $Y$ in terms of invariants of $X$ and the ramification of the map $\pi$. For the genus, we can do so using the well-known Riemann-Hurwitz formula \cite[IV.Corollary 2.4]{Hartshorne}:
$$ 2g_Y - 2 = \#G (2g_X -2) + \sum_{Q \in \pi^{-1}(S)} \sum_{i \geq 0} (\#G_i(Q) -1).$$
Here $G_i(Q)$ is the $i$-th ramification group at $Q$ with lower numbering. 

The fact that $k$ has positive characteristic gives more invariants to work with. To begin with, there is the $p$-rank $s_X$, which can be defined by the fact that $\Jac(X)[p]$ is a finite group scheme with $p^{s_X}$ geometric points. If $G$ is a $p$-group, the \emph{Deuring-Shafarevich formula} (\cite[Theorem 4.1]{SubraoDS}) expresses the $p$-rank of $Y$ in terms of the $p$-rank of $X$:
$$s_Y - 1 = \# G (s_X-1) + \sum_{Q \in \pi^{-1}(S)} (\# G_0(Q)-1).$$
In fact, by the solvability of $p$-groups, the key case is $G=\Z/p\Z$, which we focus on from here on. Such a cover is called an \emph{Artin-Schreier cover}. 

The $p$-rank can also be defined using sheaf cohomology. Denote by $\sigma$ the Frobenius morphism of $k$. Since $H^1(X,\cO_X)$ is equipped with a natural $\sigma$-linear absolute Frobenius operator $F_X$, applying Grothendieck-Serre duality gives a $\sigma^{-1}$-linear operator $V_X$ on $H^0(X,\Om_X^1)$ called the \emph{Cartier operator}. 

By a result of Oda (\cite[Corollary 5.11]{Oda}), the $p$-rank $s_X$ equals the dimension of the subspace of $H^0(X,\Om_X^1)$ on which $V_X$ acts bijectively. Since $H^0(X,\Om_X^1)$ has dimension $g_X$, the $p$-rank equals the rank of $V_X^{g_X}$ (or any higher power). Through the Deuring-Shafarevich formula the rank of $V_Y^{g_Y}$ can be expressed in terms of the rank of $V_X^{g_X}$ and the ramification invariants of $\pi$. It is natural to ask for a similar result for $V_Y^n$ with $n<g_Y$. In the case $p=2$ and $X=\PP^1$, it was shown in~\cite{Prieshyp2} that the rank of each power of $V_Y$, and in fact the entire Ekedahl-Oort type, is determined by the ramification invariants $d_Q$. However, in general the rank of $V_Y^n$ can still vary if we keep $X$ and the ramification data of $\pi:Y\to X$ constant.

Starting with the rank of $V_X$ itself gives rise to the \emph{$a$-number}, $a_X:=\dim_k \ker( V_X)$. Equivalently, $a_X$ equals $\dim_k \textup{Hom} (\alpha_p, \Jac(X)[p])$. Contrary to the genus and the $p$-rank, the $a$-number is not determined by $X$ and the ramification invariants of $\pi$, except in the special cases (see~\cite{Voloch1988},~\cite{PriesConstant}). Instead, Booher and Cais proved in \cite[Theorem 1.1]{BoCaASc} that the $a$-number lies within the following bounds. 

\begin{theorem}[Booher and Cais, 2020] \label{thmBC}
Let $\pi:Y \to X$ be an Artin-Schreier cover with branch locus $S \subset X$. For each $Q \in S$, let $d_Q$ be the unique break in the lower-numbering ramification filtration. Then for any $0 \leq j \leq p-1$ we have 
$$\sum_{Q\in S} \sum_{i=j}^{p-1} \left\lfloor \frac{id_Q}{p} \right\rfloor- \left\lfloor \frac{id_Q}{p}- \left(1-\frac{1}{p} \right) \frac{jd_Q}{p} \right\rfloor \leq a_Y \leq pa_X + \sum_{Q\in S} \sum_{i=1}^{p-1} \left\lfloor \frac{id_Q}{p} \right\rfloor- (p-i) \left\lfloor \frac{id_Q}{p^2} \right\rfloor.$$
\end{theorem}

These bounds were obtained by viewing the Cartier operator as a map of sheaves and performing computations locally. In this paper, that approach is generalized to arbitrary powers of the Cartier operator. This results in the following bounds for the \emph{$n$-th $a$-number} $$a_Y^n:=\dim \ker(V_Y^n).$$
Given integers $d$, $i$ and $n$, we define $\sigma_p(d,i,n)$ to be the number of integers $0 < l \leq \lceil id/p \rceil$ that are not $1 \bmod p^n$ and such that the first $n$ digits in the $p$-adic expansion of $(1-l)/d$ sum up to at most $p-1-i$.

\begin{theorem}[Theorem~\ref{thmmain}] \label{thmmainintro}
With the notation above, we have
\begin{equation*}
    a_Y^n \leq pa_X^n + \sum_{Q\in S} \sum_{i=1}^{p-1}\left( \left \lfloor \frac{id_Q}{p} \right\rfloor - \left\lfloor \frac{id_Q}{p^{n+1}} \right\rfloor -\sigma_p(d_Q,i,n) \right).
\end{equation*}
Moreover, for every $0\leq j \leq p-1$, we have
\begin{equation*} 
    a_Y^n \geq \sum_{Q \in S} \sum_{i=j}^{p-1} \left( \left\lfloor \frac{id_Q}{p} \right\rfloor - \left\lfloor \frac{id_Q}{p} - \left(1-\frac{1}{p^n}\right) \frac{jd_Q}{p} \right\rfloor \right).
\end{equation*}
\end{theorem}

Setting $n=1$ recovers Theorem~\ref{thmBC}. In Corollary~\ref{corsharp} it is proved that the upper bound is sharp for covers of $\PP^1$ branched at a single point. In particular, if $Y$ has the model $y^p-y=t^d$, then $a_Y^n$ attains the upper bound of Theorem~\ref{thmmainintro} for every $n$. 

When $p=2$ and $a_X^n=0$ the bounds coincide, giving the following corollary.

\begin{corollary}[Corollary~\ref{corp=2}]
If $p=2$ and $X$ is ordinary, then for $n>0$ we have
$$a_Y^n = \sum_{Q \in S} \frac{d_Q-1}{2} - \left \lfloor \frac{d_Q}{2^{n+1}} \right \rfloor.$$
\end{corollary}

In other words, the ranks of the operators $V_Y^n$ are determined by the ramification invariants. In particular, they do not depend on the location of the branch points or the lower order terms of local models at branch points.

It is a consequence of~\cite[Lemma 3.2 and 3.3]{Prieshyp2} that this corollary holds when $X$ is the projective line, but Corollary~\ref{corp=2} shows that we may replace the projective line with any ordinary curve. It remains an open question whether the entire Ekedahl-Oort type of $Y$ is determined by $X$ and the ramification invariants $d_Q$ when $p=2$ and $X$ is ordinary.

\subsection{Structure of this paper}

In Section~\ref{sec2}, we prove the required lemmas about the splitting of exact sequences involving the Cartier operator. 
In Section~\ref{sec3}, the existing theory of minimal models of Artin-Schreier covers is summarized. Any Artin-Schreier cover $\pi: Y \to X$ has a model $y^p-y=f$, where $f \in k(X)$ has a pole of order exactly $d_Q$ at each ramified point $Q \in S$. In the case $g_X>0$, it is not possible to prevent $f$ from having poles elsewhere, but it is possible to ensure that $f$ has only one other pole, at some auxiliary point $Q'$, of order $d_Q' \leq p(2g_X-2)$. A useful sheaf $\cF_0$ is constructed, which equals $\pi_* \Om_Y^1$ away from the auxiliary point $Q'$. Then explicit divisors $E_i$ are introduced such that
$$\cF_0 \cong \bigoplus_{i=0}^{p-1} \Om_X^1(E_i).$$
In Section~\ref{sec4}, the operator $V_Y^n$ is described in terms of $V_X^n$. This results in a map of $\cO_X$-modules 
$$\ph: \pi_*\ker V_Y^n \hookrightarrow \bigoplus_{i=0}^{p-1} \ker V_X^n(F_*^n(E_i+p^nD_i)).$$
Here $D_i$ are effective divisors needed to split the exact sequence obtained by applying the Cartier operator to $\Om_X^1(E_i)$. The existence of such divisors $D_i$ are guaranteed by the work in Section~\ref{sec2}. The map $\ph$ induces an isomorphism $\ph_\eta$ when localized at the generic point $\eta$ of $X$. However, on global sections it fails to be surjective in general. The remainder of the paper revolves around analyzing all the ways in which $H^0(\ph)$ can fail to be surjective. Equivalently, given an element $$\nu \in \bigoplus_{i=0}^{p-1} \ker V_X^n(F_*^n(E_i+p^nD_i)),$$ we analyze all the ways in which $\ph_\eta^{-1}(\nu)$ can fail to be regular.
In Section~\ref{sec5}, a Galois invariant filtration of sheaves is constructed for this purpose. At each step in the filtration, possible irregular behavior of $\ph_\eta^{-1}(\nu)$ is recorded. This results in Theorem~\ref{thmexseqMj}. 
In Section~\ref{sec6}, all the aforementioned tools are used to prove concrete bounds on $a_Y^n$. The main ingredients are Theorem~\ref{thmexseqMj} and the main theorem of~\cite{Tango}, which allows us to compute the dimensions of the relevant spaces of global sections. Finally, this results in the bounds of Theorem~\ref{thmmainintro}.
In Section~\ref{sec7} the bounds are applied to some examples and it is studied how the bounds behave as $n$ increases.

\subsection{Main differences with the case \texorpdfstring{$n=1$}{}}

Since the methods and results of this paper are similar in nature to those in~\cite{BoCaASc}, where the case $n=1$ is treated, an attempt has been made to stick to the structure and notation of that paper. Hopefully this will make it easier to understand for those who have already read that paper. We now outline which new ingredients are required for $n>1$ and where the added difficulty is manifested.

In Section~\ref{sec2}, Lemma~\ref{lemsecs} describes an adapted effect of the section on the order of vanishing of a differential. In Section~\ref{sec3}, the divisors $E_i$ need to have a different order at the auxiliary point $Q'$. Apart from that, Sections~\ref{sec2} and~\ref{sec3} remain largely the same as in~\cite{BoCaASc}, as the exponent $n$ of the Cartier operator is barely used. In contrast, Section~\ref{sec4} depends heavily on $n$. Therefore this section is significantly different from~\cite[Section 4]{BoCaASc}. Compare, for instance, Lemma~\ref{lemVYinVX} with~\cite[Lemma 4.1]{BoCaASc}. The difficulty introduced by this more complicated formula continues to have consequences in Section~\ref{sec5} and~\ref{sec6}. For this reason most statements need to be reproved to hold for general $n$. In particular, the skyscraper sheaves $M_j$ have a different definition (see Definitions~\ref{defMjQ} and~\ref{defM-1}) that accounts for general $n$. The technical Lemma~\ref{lem5.10} requires a significantly more complicated approach than in the case $n=1$. In Section~\ref{sec6}, the triples $(Q,l,j)$ need to be defined differently (see Definition~\ref{defQlj}), essentially due to the formula in Lemma~\ref{lemVYinVX}. This in turn accounts for the involvement of the function $\sigma_p(d_Q,i,n)$ in the upper bound for $a_Y^n$, which is perhaps not the generalization one would naively expect. As for the lower bound, Lemma~\ref{lem5.4}(i) is rather different and more complicated to prove than in the case $n=1$. For this reason the numbers $c(i,j,Q)$ need to be defined differently (see Definition~\ref{defL(X,pi)}), which ultimately accounts for the formula for the lower bound.

\section*{Acknowledgement}

The author thanks Jeremy Booher, Bryden Cais and Damiano Testa for helpful conversations. The author thanks an anonymous referee for helpful comments.

\section{Producing splittings} \label{sec2}

Let $k$ be an algebraically closed field of characteristic $p$ and let $X$ be a smooth, projective, geometrically connected curve over $k$. One views the Cartier operator $V_X$ not as a semilinear operator on $H^0(X,\Om_X^1)$, but as a map of $\cO_x$-modules $V_X: F_* \Om_X^1 \to \Om_X^1$. Similarly, for $n \geq 1$, let $V_X^n:F_*^n\Om_X^1 \to \Om_X^1$ be the $n$-th power of the Cartier operator. Using the four properties that characterize $V_X^n$ (additivity, $p^{-n}$-linearity, annihilation of exact differentials and preservation of logarithmic differentials), it is straightforward to describe the action of the Cartier operator on stalks. For $Q \in X$ and $t_Q \in \cO_{X,Q}$ a uniformizer, we have
\begin{equation} \label{eqVXnstalk}
   V_X^n\left(\sum_{i} a_i t_Q^{i-1} dt_Q\right)= \sum_{j} a_{p^nj}^{1/p^n} t_Q^{j-1} dt_Q. 
\end{equation}

Viewing $V_X^n$ as a map of sheaves gives rise to the canonical exact sequence of $\cO_X$-modules
\begin{equation} \label{exseqVXn}
    0 \to \ker V_X^n \to F_*^n \Om_X^1 \to \im V_X^n \to 0.
\end{equation}
Unfortunately, this sequence is in general not split. We first treat a special case in which it is split. 

\begin{lemma} \label{lemsplit}
If $X=\PP^1$, then \eqref{exseqVXn} is split as a sequence of $\cO_{\PP^1}$-modules.
\end{lemma}

\begin{proof}
Identify the function field of $\PP^1$ with $k(t)$. Then by elementary properties of the Cartier operator we have
$$ \ker V_X^n = \left\{ \sum_{i=1}^{p^n-1} h_i(t) t^{i-1} dt \; | \; h_i \in k(t^{p^n})  \right\}.$$
Hence \eqref{exseqVXn} has the section
\begin{align*}
    s: \im V_X^n=\Om_X^1 &\to F_*^n \Om_X^1 \\
    \sum a_i t^{i-1} dt &\mapsto \sum_i a_{i}^{p^n} t^{p^n i-1} dt.
\end{align*}
The retraction is then given by
\begin{align*}
    r: F_*^n \Om_X^1 &\to \ker V_X^n \\
    \sum a_i t^{i-1} dt &\mapsto \sum_{p^n\not |i}  a_{i} t^{i-1} dt.
\end{align*}
\end{proof}

If $X\neq \PP^1$, we need additional divisors to split \eqref{exseqVXn}. For this, work from~\cite[section 2]{BoCaASc} is still useful for our current purposes.

\begin{lemma} \label{lemexseqE}
Given any effective divisor $E=\sum_j m_j Q_j$ on $X$, we have an exact sequence
\begin{equation} \label{eqexseqE}
0 \to \ker V_X^n (F_*^n E) \to F_*^n(\Om_X^1(E)) \to \im V_X^n(\overline{E}) \to 0
\end{equation}
where $\overline{E}:= \sum_j \lceil m_j/p^n \rceil Q_j$. Each term is locally free.
\end{lemma}
\begin{proof}
This is very similar to~\cite[Lemma 2.5]{BoCaASc}. For a closed point $Q \in X$ and a germ $\om \in F_*^n(\Om_X^1(E))_Q$, we have $\ord_Q(V_X^n(\om)) \geq \lceil \ord_Q(\om)/p^n \rceil$, so the map on the right is well defined. Exactness on the left and the middle is then clear. We check exactness on the right on completed stalks using \eqref{eqVXnstalk}. Finally, since each sheaf is a subsheaf of $F_*^n\Om_{X,\eta}^1$, they are torsion-free and therefore locally free. 
\end{proof}

Using a divisor $D$, we can split the exact sequence \eqref{eqexseqE}.

\begin{lemma} \label{lemsplitting}
Let $S$ be a finite set of points on $X$ Given any effective divisor $E=\sum_j m_i Q_j$ on $X$ there exists a divisor $D=\sum_i P_i$ on $X$ with $P_i \notin S$ and a morphism 
$$ r: F_*^n(\Om_X^1(E)) \to \ker V_X^n (F_*^n E)(D)=\ker V_X^n (F_*^n(E+p^n D)$$
such that the restriction of $r$ to $\ker V_X^n(F_*^n E)$ equals the natural inclusion $$\ker V_X^n (F_*^n E) \hookrightarrow \ker V_X^n (F_*^n (E+p^n D)).$$ Equivalently, one gets a section
$$ s: \im V_X^n(\overline{E}) \to F_*^n(\Om_X^1(E+p^n D))$$
such that $V_X^n \circ s$ equals the natural inclusion $\im V_X^n(\overline{E}) \to \im V_X^n(\overline{E}+D)$.
\end{lemma}
\begin{proof}
This is a direct consequence of~\cite[Proposition 2.6]{BoCaASc}.
\end{proof}

Later we will need the following lemma about the influence of $s$ on the order of vanishing at points in the support of $E$.

\begin{lemma} \label{lemsecs}
We can choose the section
$$ s: \im V_X^n(\overline{E}) \to F_*^n(\Om_X^1(E+p^n D))$$ such that given a point $Q \in \supp (E)$ and a section $\om$ such that $\ord_Q(\om) \geq d$ we have 
$$ \ord_Q(s(\om)) \geq (d+1)p^n-1.$$
\end{lemma}
\begin{proof}
Since $V_X^n \circ s$ is the natural inclusion, we know that $V_X^n(s(w))=w$. Now, let $t_Q$ be a uniformiser of $\cO_{X,Q}$. Suppose that $\ord_Q(s(\om)) < (d+1)p^n-1.$ We may pick $s$ such that $s(\om)$ never has terms that lie in the kernel of $V_X^n$. This means that $s(\om)$ has a term of the form $t_Q^{p^nl-1} dt_Q$, with $p^nl-1<(d+1)p^n-1$, meaning $l<d+1$. Then \eqref{eqVXnstalk} implies that $\om=V_X^n(s(\om))$ has a term $t_Q^{l-1}$, which contradicts the assumption $\ord_Q(\om)\geq d$.
\end{proof}

\section{Artin-Schreier covers and differential forms} \label{sec3}

This section is concerned with finding a model for our Artin-Schreier cover $\pi: Y \to X$. Since this investigation does not involve the Cartier operator, we can make use of the results in~\cite[section 3]{BoCaASc}, with a slight adjustment. We give a brief summary of these results, as they will be important later. Recall that the cover has branch locus $S \subset X$ and for each $Q \in S$, $d_Q$ is the unique break in the lower-numbering ramification filtration. 

In the case $X=\PP^1$, every Artin-Schreier extension of $k(X)$ has a model of the simplest possible form:
$$ Y: y^p-y=f,$$
where $f$ has exactly the prescribed poles: at any $Q\in S$, it has a pole of degree $d_Q$. When $X$ is not the projective line, it is not possible to prevent $f$ from having poles outside $S$, but it is possible to force $f$ to have just one pole outside $S$. This is the auxiliary point $Q'$, where $f$ has a pole of order $d_{Q'} \leq p(2g_X-2)$, with $p | d_{Q'}$. The existence of such an element $f\in k(X)$ is proved in~\cite[Lemma 3.2]{BoCaASc}. We fix such an $f$ for the rest of this paper.

Using this model for $Y$, we construct a Galois filtration of sheaves on $X$. We need the following definitions:

\begin{align*}
    S'&:=S \cup \{Q'\} \\
    d_Q&:= \min \{0, \ord_f(Q) \} \\
    n_{Q,i} &:= \begin{cases}
    \left \lceil \frac{d_Q(p-1-i)}{p} \right \rceil & \hbox{if $Q \in S$} \\
    p^{n-1}d_{Q'}(p-1-i) & \hbox{if $Q = Q'$} \\
    0 & \hbox{else}
    \end{cases} \\
    E_i &:= \sum_{Q \in S'} n_{Q,i} Q \\
    \overline{E}_i &:= \sum_{Q\in S} \lceil n_{Q,i}/p \rceil Q. \\
\end{align*}

These definitions are all the same as in~\cite{BoCaASc}, except for $n_{Q',i}$. In the case $n=1$, the definitions coincide. In Remark~\ref{rmknQ'i} it will be explained why this generalization of $n_{Q',i}$ is the right one for our purposes. 

Because of the auxiliary point $Q'$, the Galois filtration on $\pi_* \Om_Y^1$ will not have the desired properties. Instead we must allow poles at $Q'$. To be precise, we work with a sheaf $\cF_0 \subset (\pi_*\Om_Y^1)_\eta$, where $\eta$ is the generic point of $X$, whose global sections on an open $U \subseteq X$ are
$$\cF_0(U) = \left\{ \om=\sum_i \om_i y^i \, : \, \ord_Q(\om_i) \geq -n_{Q,i} \text{ for $Q \in U$ and $0 \leq i \leq p-1$}\right\}.$$
The condition $\ord_Q(\om_i) \geq -n_{Q,i}$ ensures that $\om$ is regular away from $Q'$. At the point $Q'$, we allow a pole of order $n_{Q',i}=p^{n-1}d_{Q'}(p-1-i)$. In the case $X=\PP^1$, the point $Q'$ is not needed, so that we have $\cF_0 = \pi_*\Om_Y^1$. 

We now fix a generator $\tau$ of $\Gal(k(Y)/k(X))$ and define the Galois filtration on $\cF_0$ as $$W_i= \ker\left ((\tau-1)^{i+1}: \cF_0 \to \cF_0 \right).$$ Note that $W_{-1}=0$ and $W_{p-1}=\cF_0$. One then obtains split exact sequences
\begin{equation} \label{exseqgal}
0 \to W_{i-1} \to W_i \to \Om_X^1(E_i) \to 0.
\end{equation}
for each $0 \leq i \leq p-1$. This fact is proved in~\cite[Proposition 3.12]{BoCaASc}. Although our divisors $E_i$ have a different order at $Q'$, the proof is not affected by this.
As a result of these split exact sequences, we obtain an isomorphism
\begin{equation} \label{eqsumgal}
\cF_0 \cong \bigoplus_{i=0}^{p-1} \Om_X^1(E_i).
\end{equation}
This isomorphism is the content of~\cite[Cor. 3.14]{BoCaASc}.

\section{Powers of the Cartier operator on stalks} \label{sec4}

From now on we assume, using the results summarized in the previous section, that the Artin-Schreier cover $\pi: Y \to X$ has a minimal model $y^p-y=f$, where $f \in k(X)$ has a pole of order $d_Q$ at each $Q\in S$ and a pole of order $d_{Q'}$ at an auxiliary point $Q'\in X \setminus S$. In this section we analyze the relation between $V_Y^n: F_*^n \Om_Y^1 \to \Om_Y^1$ and $V_X^n: F_*^n \Om_X^1 \to \Om_X^1$ on stalks at the generic point. In particular, Proposition~\ref{propphieta} establishes an isomorphism
$$ (\ker V_Y^n)_{\eta'} \cong (\pi_*\ker V_Y^n)_\eta \cong \bigoplus_{i=0}^{p-1} (\ker V_X^n)_\eta,$$
where $\eta'$ is the generic point of $Y$. This isomorphism gives rise to an inclusion of sheaves $$\ph: \pi_*\ker V_Y^n \hookrightarrow \bigoplus_{i=0}^{p-1} \ker V_X^n(F_*^n(E_i+p^nD_i)). $$ Analyzing the image of this map carefully eventually results in the bounds of Theorem~\ref{thmmainintro}.

We repeatedly use the following four properties that characterize the Cartier operator:
\begin{enumerate}
    \item (Additivity) $V_X(\om_1+\om_2)=V_X(\om_1)+V_X(\om_2)$ for $\om_1,\om_2 \in H^0(X,\Om_X^1)$.
    \item ($p^{-1}$-linearity) $V_X(f^p\om)=fV_X(\om)$ for $f\in k(X)$ and $\om\in H^0(X,\Om_X^1)$.
    \item (Annihilation of exact differentials) $V_X(df)=0$ for $f \in k(X)$.
    \item (Preservation of logarithmic differentials) $V_X(df/f)=df/f$ for $f \in k(X)$.
\end{enumerate}

Given $l,m \in \Z_{\geq 0}$, we define the operator $W_{l,m}$ on $H^0(X,\Om_X^1)$ by
$$W_{l,m}(\om) :=  V_X\left( \binom{l}{m} (-f)^{l-m} \om \right).$$

\begin{lemma} \label{lemVYinVX}
Let $\om \in \Om_{Y,\eta'}^1$. Write $\om=\sum_{i=0}^{p-1}\om_i y^i$ with $\om_i \in \Om_{X,\eta}^1$. We have 
$$V_Y^n(\om) = \sum_{j=0}^{p-1} \left(
\sum_{j\leq j_{1}\leq \ldots \leq j_n \leq p-1}
W_{j_{1},j}( W_{j_{2},j_{1}}( \ldots W_{j_{n},j_{n-1}}(\om_{j_n}) \ldots )) \right) y^j$$
\end{lemma}
\begin{proof}
We first determine the action of $V_Y^n$ on the summands $\om_i y^i$, so that we can deduce the result by additivity. We claim
$$V_Y^n(\om_i y^i) =
\sum_{0 \leq j_1 \leq \ldots \leq j_n \leq i}
W_{j_{2},j_1}( W_{j_{3},j_{2}}( \ldots W_{i,j_n}(\om_{i}) \ldots ))y^{j_1}.$$
One observes that the case $n=0$ is trivial and proceeds by induction. 
\begin{align*}
V_Y^n(\om_i y^i) &= V_Y^n(\om_i (y^p-f)^i) \\
&= V_Y^n\left( \sum_{j_n=0}^{i} \binom{i}{j_n} \om_i y^{pj_n} (-f)^{i-j_{n}} \right) \\
&= \sum_{j_n=0}^{i} V_Y^{n-1} \left( V_Y \left( \binom{i}{j_n} \om_i y^{pj_n} (-f)^{i-j_n}\right) \right) \\
&= \sum_{j_n=0}^{i} V_Y^{n-1} \left( V_X \left(\binom{i}{j_n} (-f)^{i-j_n} \om_i\right)y^{j_n}\right) \\
&= \sum_{j_n=0}^i V_Y^{n-1} \left( W_{i,j_n} (\om_i) y^{j_n} \right) \\
&= \sum_{j_n=0}^i \sum_{0 \leq j_1 \leq \ldots \leq j_{n-1} \leq j_n}
W_{j_{2},j_1}( W_{j_{3},j_{2}}( \ldots W_{i,j_n}(\om_{i}) \ldots ))y^{j_1}.
\end{align*}
Then the result follows from collecting $y^j$ terms.
\end{proof}

This explicit description of the action of the iterated Cartier operator allows us to describe its kernel.

\begin{corollary} \label{corVYker}
Let $\om = \sum_{i=0}^{p-1} \om_i y^i$ be an element of $\Om_{Y,\eta'}^1$. Then $V_Y^n(\om)=0$ if and only if $V_X^n(\om_{p-1})=0$ and for every $j<p-1$ we have
$$ V_X^n(\om_j) = - \sum_{i=j+1}^{p-1} \left( \sum_{j \leq j_1 \leq \ldots \leq j_{n-1} \leq i} W_{j_1,j}(W_{j_2,j_1}(\ldots W_{i,j_{n-1}}(\om_i)\ldots))\right) . $$
\end{corollary}
\begin{proof}
This is an application of Lemma~\ref{lemVYinVX}. By assumption each $y^j$ term of $V_Y^n(\om)$ vanishes. If $j=p-1$, the coefficient equals $V_X^n(\om_{p-1})$, which must therefore equal zero. In the case $j<p-1$ we take the term with $j=j_1=j_2=\ldots=j_n$, which equals $V_X^n(\om_j)$, out of the large sum, giving the result. 
\end{proof}

Recall from Lemma~\ref{lemsplitting} that we have the maps
\begin{align*}
    s_i: \im V_X^n(\overline{E}_i) &\to F_*^n \Om_X^1(E_i+p^nD_i) \\
    r_i: F_*^n \Om_X^1(E_i) &\to \ker V_X^n(F_*^n(E_i+p^nD_i)).
\end{align*}
Lemma~\ref{lemsplit} shows that we can take $D_i=0$ in the case $X=\PP^1$. Using these maps, we can write any element $\om_i \in F_*\Om_X^1(E_i)$ as
$$ \om_i= r_i(\om_i) + s_i(V_X^n(\om_i)).$$

\begin{definition} \label{defphi}
We define the map $\ph$ of sheaves on $X$:
$$\ph: \pi_*\ker V_Y^n \hookrightarrow F_*^n \cF_0 \cong \bigoplus_{i=0}^{p-1} F_*^n(\Om_X^1(E_i)) \to \bigoplus_{i=0}^{p-1} \ker V_X^n(F_*^n(E_i+p^nD_i)),$$
where the rightmost arrow is $\bigoplus r_i$.
\end{definition}

\begin{proposition} \label{propphieta}
On the stalk at the generic point $\eta$ of $X$ $$\ph_\eta: \left( \pi_*\ker V_Y\right)_{\eta} \to \bigoplus_{i=0}^{p-1} \left(\ker V_X^n(F_*^n(E_i+p^nD_i))\right)_\eta= \bigoplus_{i=0}^{p-1} \left(\ker V_X^n\right)_\eta$$ is an isomorphism of $K(X)$-vector spaces.
\end{proposition}
\begin{proof}
Using Corollary~\ref{corVYker}, we explicitly construct the inverse of $\ph_\eta$. Let $$\nu=(\nu_0, \nu_1, \ldots , \nu_{p-1})\in \bigoplus_{i=0}^{p-1} \left(\ker V_X^n \right)_\eta$$ be arbitrary. We construct a differential $\om=\sum_{j=0}^{p-1} \om_i y^j \in \Om_Y^1$ as follows. We begin by setting $\om_{p-1}=\nu_{p-1}$ and then define the remaining $\om_j$'s by descending induction:
\begin{equation} \label{eqphiinv}
\om_j= \nu_j + s_j\left(  - \sum_{i=j+1}^{p-1} \left( \sum_{j \leq j_1 \leq \ldots \leq j_{n-1} \leq i} W_{j_1,j}(\ldots W_{i,j_{n-1}}(\om_i) \ldots )\right) \right).
\end{equation}
Then Corollary~\ref{corVYker} implies $\om \in \left( \pi_*\ker V_Y\right)_\eta.$ The assignment $\nu \mapsto \om$ is the desired inverse of $\ph_\eta$.
\end{proof}

We now analyze an example where a non-leading term of $f$ increases the rank of the square of the Cartier operator. Upon closer inspection, the non-leading term accounts for an element of $\bigoplus_{i=0}^{p-1} H^0(X,\ker V_{X}^2(F_*^2 E_i))$ that fails to be in the image of $H^0(\ph)$. 

\begin{example} \label{eg7}
We present two covers of $X=\PP^1$ in characteristic $7$ with the same ramification data on which the square of the Cartier operator has a different rank. Define $f_1:= t^{-4}$ and $f_2:=t^{-4}+t^{-3}$. For $i=1,2$, let $Y_i$ be defined by $y^p-y=f_i$. Both curves have genus $9$, as determined by the Riemann-Hurwitz formula. On $Y_1$, the Cartier operator is zero, so $a_{Y_1}^n=9$ for each $n$. On the other hand, $Y_2$ has $a$-number $6$. Furthermore, $V_{Y_2}^2$ still has rank $1$. Thus we have $a_{Y_1}^2=9$ and $a_{Y_2}^2=8$. Let's investigate this discrepancy more closely with the amassed tools.

We have but one branch point $Q=0$, with $d_Q=4$. We compute $n_{Q,i}=\lceil 4(6-i)/7 \rceil$ for $0 \leq i \leq p-1$.
Since $t^{-1}dt$ is preserved by the Cartier operator $V_{\PP^1}$, and therefore also by all its powers, one has
$$H^0(\PP^1,\ker V_{\PP^1}^n(F_*^n E_i))= \{t^{-l} dt \; | \; 2 \leq l \leq n_{Q,i} \}.$$
Here it is used that $n_{Q,i} < p+1$ for every $i$.
Let $v_{j,l} \in \bigoplus_{i=0}^{p-1} H^0(\PP^1,\ker V_{\PP^1}^n(F_*^n E_i))$, with $2\leq l \leq n_{Q,j}$, denote the element that has $t^{-l}dt$ in the $j$-th position and $0$ elsewhere. Then $\bigoplus_{i=0}^{p-1} H^0(\PP^1,\ker V_{\PP^1}^n(F_*^n E_i))$ is spanned by 
$$\{\nu_{0,2},\nu_{0,3}, \nu_{0,3}, \nu_{1,2}, \nu_{1,3}, \nu_{2,2}, \nu_{2,3}, \nu_{3,2}, \nu_{4,2} \}.$$

We now specialize to $n=2$. Proposition~\ref{propphieta} provides that the map on global sections
$$ H^0(\ph): H^0(\PP^1,\pi_* \ker V_Y^2) \to \bigoplus_{i=0}^{p-1} H^0(\PP^1,\ker V_{\PP^1}^2(F_*^2 E_i))$$
is injective. We show that $H^0(\ph)$ fails to be surjective if $Y=Y_2$. Since $n_{Q,4}=2$ we infer that $\bigoplus_{i=0}^{p-1} H^0(\PP^1,\ker V_{\PP^1}^2(F_*^2 E_i))$ contains the element 
$$\nu_{4,2}:= (0, 0, 0, 0, t^{-2}dt, 0, 0).$$
Viewing $\nu_{4,2}$ as an element of $\bigoplus_{i=0}^{p-1} \left(\ker V_{\PP^1}^2(F_*^2(E_i))\right)_\eta$, we utilize Proposition~\ref{propphieta} to compute $\om=\sum_{j=0}^{p-1}\om_j y^j:=\ph_\eta^{-1}(\nu_{4,2})$. The formula \eqref{eqphiinv} in the present situation reads
$$\om_j = (\nu_{4,2})_j + s_j\left( - \sum_{i=j+1}^{p-1} \left( \sum_{j \leq j_1 \leq i} V_{\PP^1}\left(\binom{j_1}{j} (-f)^{j_1-j} V_{\PP^1}\left(\binom{i}{j_1} (-f)^{i-j_1} \om_i\right)\right)   \right)  \right)$$
It is a straightforward application that $\om_j=0$ for $j>4$ and $\om_4=t^{-2}dt$. Now, for $j<4$, we iterate the use of this formula to show that $\om_j=0$ for $0<j<4$: there are no values of $j$ and $j_1$ that enable $\om_i$ to survive two iterations of the Cartier operator. That changes in the computation of $\om_0$. In the case $Y=Y_1$, $(-f)$ only contributes a term $t^{-4}$. In the case $Y=Y_2$, one has
\begin{align*}
    V_{\PP^1}((-f_2)^2 t^{-2}dt) &= V_{\PP^1}((t^{-8}+2t^{-7}+t^{-6})t^{-2}dt) \\
    &= t^{-1}\left(V_{\PP^1}(t^{-3}dt)+2V_{\PP^1}(t^{-2}dt)+V_{\PP^1}(t^{-1}dt) \right) \\
    &= t^{-2}dt.
\end{align*} 
Hence $(i,j_1)=(4,2)$ yields the summand
\begin{align*}
    V_{\PP^1}\left(\binom{2}{0} (-t^{-4}+t^{-3})^{2} V_{\PP^1}\left(\binom{4}{2} (-t^{-4}+t^{-3})^{2} t^{-2}dt\right)\right)&= V_{\PP^1}\left(6(t^{-4}+t^{-3})^2t^{-2}dt \right) \\
    &=6t^{-2}dt.
\end{align*}
One verifies that all other summands are zero and hence
$\ph_\eta^{-1}(\nu_{4,2})=t^{-2}dt y^4 + 6t^{-50}dt$, which fails to be regular at $Q=0$. Hence $H^0(\ph)$ is not a surjection and therefore 
$$\dim \ker V_{Y_2}^2 < 9.$$
A straightfoward computation shows that this phenomenon does not happen to basis vectors other than $\nu_{4,2}$. We thus conclude
\begin{align*}
    a_{Y_1}^2&=\dim \ker V_{Y_1}^2 = 9 \\
    a_{Y_2}^2&=\dim \ker V_{Y_2}^2 = 8.
\end{align*}
Since $V_{Y_2}$ is nilpotent, we have $a_{Y_2}^3=9$.
\end{example}

\begin{remark}
An example can also be found in characteristic $5$, with
\begin{align*}
    f_1&:=t^{-6} \\
    f_2&:=t^{-6}+t^{-4}.
\end{align*}
We then have $g_1=g_2=a_1=10$, but $V_{Y_2}^2$ has rank $2$. Another example is given in characteristic $3$:
\begin{align*}
    f_1&:=t^{-10} \\
    f_2&:=t^{-10}+t^{-8}.
\end{align*}
The resulting curves also have genus $9$, but $V_{Y_2}^2$ again has rank $2$.
\end{remark}

\section{Short Exact Sequences and the Kernel} \label{sec5}

In the previous section, we constructed a map of sheaves
$$\ph: \pi_*\ker V_Y^n \hookrightarrow F_*^n \cF_0 \cong \bigoplus_{i=0}^{p-1} F_*^n(\Om_X^1(E_i)) \to \bigoplus_{i=0}^{p-1} \ker V_X^n(F_*^n(E_i+p^nD_i)).$$
This map induces an isomorphism $\ph_\eta$ at the generic point $\eta$, which is the content of Proposition~\ref{propphieta}. Therefore the induced map on global sections
$$ H^0(\ph): H^0(X, \pi_* \ker V_Y^n) \to \bigoplus_{i=0}^{p-1} H^0\left(X,\ker V_X^n(F_*^n(E_i+p^nD_i))\right)$$
is injective. However, in general it fails to be surjective, as is illustrated in Example~\ref{eg7}. Since we are interested in $a_Y^n = \dim_k H^0 (X, \pi_* \ker V_Y^n)$, our next step is to measure all the ways in which $H^0(\ph)$ can fail to be surjective. Equivalently, given an element $\nu$ of $\bigoplus_{i=0}^{p-1} H^0(X,\ker V_X^n(F_*^n(E_i+p^nD_i)))$, we measure how $\ph_\eta^{-1}(\nu)$ can fail to be regular.

In this section we introduce the tools necessary to do this. We construct a filtration of sheaves
\begin{equation} \label{eqGfilt}
0 \subset \im \ph = \G_{-1} \subset \G_0 \subset \ldots \subset \G_p = \bigoplus_{i=0}^{p-1} \ker V_{X^n}(F_*^n(E_i+p^nD_i).
\end{equation}
The rightmost sheaf and its global sections will be accessible to us (see Proposition~\ref{proptango}), while the leftmost sheaf $\G_{-1}=\im \ph$ is what we're interested in. The goal of this section is to work our way from right to left, by computing the associated graded sheaves of the filtration \eqref{eqGfilt}, which will be skyscraper sheaves. This is done through local computations at the branch points of the cover, the support of the auxiliary divisors $D_j$ required to split \eqref{exseqVXn} and the auxiliary point $Q'$.

\subsection{Filtration and Skyscraper Sheaves}

We focus first on defining the sheaves $\G_j$ in \eqref{eqGfilt}.
For an abelian group $A$ and a point $P\in X$, let $P_*A$ denote the sheaf on $X$ whose stalk is $A$ at $P$ and $0$ away from $P$. Similarly, given a group homomorphism $h:A\to B$, we let $P_*h:P_*A \to P_*B$ be the map of sheaves on $X$ that is $h$ at $P$ and trivial away from $P$. 

\begin{definition} \label{defggj}
We let $S_j'=S' \cup \supp (D_j)$. Let $r_{Q,j}:=n_{Q,j}$ if $Q \in S'$ and $r_{Q,j}:=p^n$ if $Q \in \supp (D_j)$. We define 
$$ \iota: \bigoplus_{i=0}^{p-1} \ker V_X^n(F_*^n(E_i+p^nD_i)) \hookrightarrow \bigoplus_{i=0}^{p-1} (\ker V_X^n)_\eta \cong (\ker V_Y^n)_{\eta'} \hookrightarrow F_*^n(\pi_*\Om_Y^1)_\eta \cong \bigoplus_{i=0}^{p-1}F_*^n\Om_{X,\eta}^1.$$
$$ \rho_j: \bigoplus_{i=0}^{p-1}F_*^n \Om^1_{X,\eta} \to F_*^n\Om^1_{X,\eta} \to F_*^n \bigoplus_{Q \in S_j'} Q_* \left( \Om^1_{X,Q} [1/t_Q] /t_Q^{-r_{Q,j}} \Om^1_{X,Q} \right) . $$
Here the isomorphism $\bigoplus_{i=0}^{p-1} (\ker V_X^n)_\eta \cong (\ker V_Y^n)_{\eta'}$ is $\ph_\eta^{-1}$ from Proposition~\ref{propphieta}. The map $\bigoplus_{i=0}^{p-1}F_*^n \Om^1_{X,\eta} \to F_*^n\Om^1_{X,\eta}$ projects to the $j$-th coordinate. Let $g_j:= \rho_j \circ \iota$.
\end{definition}

Using the maps $g_j$, we now define the subsheaves $\G_j$ in \eqref{eqGfilt}.

\begin{definition} \label{defGj}
Let $0 \leq j \leq p$. For $\nu=(\nu_0, \ldots, \nu_{p-1}) \in \bigoplus_{i=0}^{p-1} \ker V_X^n(F_*^n(E_i+p^nD_i))$, let $\iota(\nu)=(\om_0, \ldots , \om_{p-1})$. Define the $\cO_X$-module
$$\G_j \subset \bigoplus_{i=0}^{p-1} \ker V_X^n(F_*^n(E_i+p^nD_i)) $$ that has sections $\nu$ with the property that $\om_i$ is a section of $F_*^n(\Om_X^1(E_i))$ whenever $j\leq i$. Define $\G_{-1}=\im \ph$.
\end{definition}

Note that $\G_j$ imposes stricter conditions on the $\om_i$ for $i \leq j$, in terms of regularity. In this way the possible failure of $\om$ to be regular can be measured step by step.   

\begin{lemma} \label{lemGjbasics}
Let $0 \leq j \leq p-1$. Recall the definition $g_j=\rho_j \circ \iota$ in Definition~\ref{defggj}.
\begin{enumerate}[(i)]
    \item $\G_j=\ker(g_{j+1} |_{\G_{j+1}})$.
    \item For $Q \notin S_j'$, we have an isomorphism $\G_{j,Q} \cong \G_{j+1,Q}$.
    \item The action of $G=\Z/p\Z$ on $(F_*^n \pi_*\Om_Y^1)_\eta \cong \bigoplus_i (\ker V_X^n)_\eta$ preserves $\G_j$.
    \item $\G_0 \cong \pi_* \ker(V_Y^n |_{\cF_0})$.
    \item $\G_{-1} \cong \pi_* \ker V_Y^n$.
    \item $\G_p=\bigoplus_{i=0}^{p-1} \ker V_X^n(F_*^n(E_i+p^nD_i))$.
\end{enumerate}
\end{lemma}

\begin{proof}
The proof of this lemma remains the same as in~\cite[Lemma 5.3]{BoCaASc}, as the exponent of the Cartier operator is never used. 
\end{proof}

Note that $\G_{-1}=\G_0$ in the case $X=\PP^1$.

\begin{lemma} \label{lem5.4}
Let $\nu:=(\nu_0, \ldots , \nu_{p-1}) \in H^0(X,\G_p)$ and set $(\om_0, \ldots , \om_{p-1}):= \iota(v)$. Fix a value $0 \leq j \leq p-1$.
\begin{enumerate}[(i)]
    \item Fix a $Q \in S$ and suppose $\nu_i=0$ for $i>j$ and $\ord_Q(\nu_j) \geq -l$ for some integer $l$. Let $\mu_{Q,i}$ be the largest multiple of $p^n$ such that $\mu_{Q,i}+1 \leq l+p^{n-1}d_Q(j-i)$. Then $\om_i=0$ for $i>j$, $\ord_Q(\om_j) \geq -l$ and for $i<j$ we have
    $$\ord_Q(\om_i) \geq \min\{\ord_Q(\nu_i), -\mu_{Q,i}-1\} \geq -l-p^{n-1}d_Q(j-i). $$
    \item For $Q \in \supp(D_j)$, we have $\ord_Q(\om_j) \geq -p^n$.
    \item $\ord_{Q'}(\om_j) \geq -p^{n-1}d_{Q'}(p-1-j)$.
\end{enumerate}
\end{lemma}

\begin{proof}
For (i), we note that $\nu_i=\om_i=0$ for $i>j$ by the construction of $\ph_\eta^{-1}$ in the proof of Proposition~\ref{propphieta}. We also have $\nu_j=\om_j$, so in particular $\ord_Q(\nu_j)=\ord_Q(\om_j) \geq -l$. For $i<j$ we proceed by decreasing induction on $i$. We obtain
\begin{align*}
    \ord _Q(\om_i)&=\ord _Q\left(\nu_i+  s_i\left(  - \sum_{m=i+1}^{p-1} \left( \sum_{i \leq j_1 \leq \ldots \leq j_{n-1} \leq m} W_{j_1,i}(\ldots W_{m,j_{n-1}}(\om_m)\ldots )\right) \right) \right) \\
    \geq \min &\left\{\ord _Q(\nu_i), \min_{ i \leq j_1 \leq \ldots \leq j_n \leq m >i} \left\{ \ord_Q  \left( s_i\left(   W_{j_1,i}(\ldots W_{m,j_{n-1}}(\om_m)\ldots )\right) \right) \right\} \right\}.
\end{align*}
By definition of $\G_p$ and $E_i$, we obtain 
$$\ord_Q(\nu_i) \geq -\left\lceil \frac{(p-1-i)d_Q}{p} \right\rceil \geq -d_Q.$$ 
We now compute the latter order of vanishing using the induction hypothesis, which provides $\ord_Q(\om_m) \geq -l-p^{n-1}d_Q(j-m)$. Note the minimal order of vanishing (i.e. the highest pole order) is obtained if we apply the Cartier operator first and then multiply by $(-f)$, since iterations of the Cartier operator roughly divide the order by $p$ (see \eqref{eqVXnstalk}). In Lemma~\ref{lemVYinVX} this amounts to the term with $j_1=m$. For any $m$, we obtain
$$\ord_Q\left( V_X\left( \binom{m}{i} (-f)^{m-i} V_X^{n-1}(\om_m)  \right) \right) \geq \left( -d_Q(m-i) + \frac{\ord_Q(\om_m)+1}{p^{n-1}}  \right)/p -1$$
\begin{align*}
    &\geq \left( -d_Q(m-i) + \frac{1-l-p^{n-1}d_Q(j-m)}{p^{n-1}}  \right)/p  -1 \\
    &=\left( (1-l)/p^{n-1} - d_Q(j-i) \right)/p-1.
\end{align*}
Then Lemma~\ref{lemsecs} provides
\begin{align*}
    \ord_Q\left(s_i \left( V_X\left( \binom{m}{i} (-f)^{m-i} V_X^{n-1}(\om_m)  \right) \right) \right) &\geq -\mu_{Q,i}-1 \\
    &\geq -l-p^{n-1}d_Q(j-i)
\end{align*}
as desired.

For (ii), note that $D_j$ is the sum of distinct points that are not in the support of $E_j$. Therefore we have $\ord_Q(\nu_j) \geq -p^n$ and $\ord_Q(f)=0$. Then the result follows from descending induction on $j$.

Similarly we prove (iii) by descending induction on $j$. By definition of $E_j$, we have $$\ord_{Q'}(\nu_j) \geq -p^{n-1}d_{Q'}(p-1-j).$$ Moreover, we have $\ord_{Q'}(f)=-d_{Q'}$, so we obtain
$$ \ord_{Q'}\left(s_j \left( V_X\left( \binom{m}{j} (-f)^{m-j} V_X^{n-1}(\om_m)  \right) \right) \right) \geq -p^{n-1}d_{Q'}(m-j) -p^{n-1}d_{Q'}(p-1-m),$$
as desired.
\end{proof}

\begin{corollary} \label{cor5.5}
Let $\nu=(\nu_0,\ldots ,\nu_{p-1}) \in H^0(X,\G_p)$ and set $(\om_0, \ldots ,\om_{p-1}):= \iota(\nu)$. For $Q\in S$ we have $$ \ord_Q(\om_i) \geq -p^{n-1}d_Q(p-1-i).$$
\end{corollary}
\begin{proof}
Note that by definition $\nu_{p-1}$ is regular. Then take $j=p-1$ and $l=0$ in Lemma~\ref{lem5.4}.(i).
\end{proof}

\begin{remark} \label{rmknQ'i}
Lemma~\ref{lem5.4}.(iii) is the reason we had to alter the definition of $n_{Q',i}$, and thereby $E_i$. It is desirable to do the analysis involving $Q'$ all in one step of the filtration, namely the step between $\G_{-1}$ and $\G_0$. We do so by allowing poles of a certain order at $Q'$ in the rest of the filtration. Lemma~\ref{lem5.4}.(iii) tells us what these allowed pole orders are. This is used in Lemma~\ref{lemMjQ'}.
\end{remark}

\begin{definition} \label{defMjQ}
    For $0 \leq j \leq p-1$ and $Q \in S$, define
    \begin{align*}
        m_{Q,j}&:=p^{n} \left\lfloor \frac{d_Q(p-1-j)}{p} \right\rfloor=p^n(n_{Q,j}-1) \\
        \beta_{Q,j}&:=\lfloor (m_{Q,j}-n_{Q,j})/p^n \rfloor+1= \lfloor (p^n-1)n_{Q,j}/p^n \rfloor. 
    \end{align*}
    Then define the $\cO_{X,Q}$-modules
    \begin{equation} \label{eqMjQ}
        M_{j,Q} := \begin{cases}  
        k[\![t_Q^{p^n}]\!]/(t_Q^{p^n\beta_{Q,j}}) & \mbox{if $Q \in S$}  \\
        k[\![t_Q^{p^n}]\!]/(t_Q^{p^n(p^n-1)})  & \mbox{if $Q \in \supp(D_j)$} \\
        0  &\mbox{if $Q=Q'$.}
        \end{cases}
    \end{equation}
    Then we define the skyscraper sheaf $M_j$ as $$M_j:= \bigoplus_{Q \in S_j'} Q_* (M_{j,Q}).$$ 
\end{definition}

We will show that for each $0 \leq j \leq p-1$ we have an exact sequence 
$$ 0 \to \G_{j} \to \G_{j+1} \to M_j \to 0.$$
We prove this by analyzing the stalks at ramified points and points in the support of $D_j$. Finally, we compute the cokernel of the inclusion $\G_{-1} \hookrightarrow \G_0$, which only involves calculations at the point $Q'$. This together gives Theorem~\ref{thmexseqMj}.

\begin{remark}
Note that $m_{Q,j}$ is the largest multiple of $p^n$ that does not exceed the integer $p^{n-1}(p-1-j)d_Q-1$ and $\beta_{Q,j}$ is the number of multiples of $p^n$ between $-m_{Q,j}$ and $-n_{Q,j}$. Let $\nu$ be a section of $\G_{j+1}$, with $(\om_0, \ldots , \om_{p-1}) :=\iota(\nu)$. Then the skyscraper sheaf $M_j$ records the coefficients in $\om_j$ of the poles with order exceeding $n_{Q,i}$. In other words, it records the ways in which $\om_j$ could fail to be a section of $\Om_X^1(E_j)$.
\end{remark}

\subsection{Local calculations above ramified points} \label{secbrp}

The goal of this section is to show that the map $\G_{j+1} \to M_{j}$ is surjective when localized at a point where $\pi:Y\to X$ ramifies.

We fix $0 \leq j \leq p-1$ and let $Q \in S$. Let $t_Q$ be a uniformizer of $\cO_{X,Q}$. We consider $\nu:=(\nu_0,\ldots,\nu_{p-1}) \in \G_{p,Q}$ and put $(\om_0, \ldots, \om_{p-1}):=\iota(\nu)$. We define $\om:=\sum_{i=0}^{p-1} \om_i y^i$. Then we have $\ph(\om)=\nu$ and $V_Y^n(\om)=0$. Moreover, by \eqref{eqphiinv} we have $\om_j=\nu_j+s_j(\xi)$, with
$$\xi := - \sum_{m=j+1}^{p-1} \left( \sum_{j \leq j_1 \leq \ldots \leq j_{n-1} \leq m} W_{j_2,i}(\ldots W_{m,j_{n-1}}(\om_m)\ldots )\right).$$

By Corollary~\ref{cor5.5} we obtain the expansion
$$ s_j(\xi) = \sum_{l=-m_{Q,j}/p^n}^{\infty} a_{l} s_j \left( t_Q^l \frac{dt_Q}{t_Q} \right) $$
in the completion $\cO^\wedge_{X,Q}$. 

\begin{definition} \label{defg'brp}
Analogous to~\cite[Definition 5.8]{BoCaASc}, we define $c_{j,Q}: \im (g_{j+1})_Q \to M_{j,Q}$ as follows:

$$ c_{j,Q}: \sum_{-m_{j,Q}/p^n}^{\infty} a_{l} s_j\left(t_Q^l \frac{dt_Q}{t_Q} \right) \bmod t_Q^{-n_{Q,j}} k[\![t_Q]\!]dt_Q \mapsto \sum_{l=0}^{\beta_{Q,j}-1} a_{-m_{Q,j}+l} t_Q^{lp^n} \bmod t_Q^{p^n\beta_{Q,j}} k[\![t_Q^{p^n}]\!]. $$

We then define the composite map 
\begin{equation} \label{eqgj+1def}
\begin{tikzcd}
{g'_{j+1,Q}:\G_{j+1,Q}} \arrow[r, "(g_{j+1}|_{\G_{j+1}})_Q"] & {\im (g_{j+1}|_{\G_{j+1}})_Q} \arrow[r, "{c_{j,Q}}"] & {M_{j,Q}}.
\end{tikzcd}
\end{equation}
\end{definition}

To prove that $g'_{j+1,Q}$ is surjective, we need the following technical result. This result and its proof are similar in nature to~\cite[Lemma 5.10]{BoCaASc}, but instead of defining each term of $\om_{p-1}$, we define only a sparse subset of these coefficients. A rather different condition of the pole order of $\om_i$ is needed, essentially because of Lemma~\ref{lemVYinVX}. This new condition is also the one that makes the adapted version of Proposition~\ref{propbrpsurj} work.

\begin{lemma} \label{lem5.10}
Suppose $\om_0, \ldots , \om_{p-2} \in k((t_Q))dt_Q$ satisfy $$\ord_Q(\om_i) \geq -p^{n-1}((p-1-i)d_Q-1)-1$$ for $0 \leq i \leq p-2$. Then there exists an $\om_{p-1} \in k[\![t_Q]\!]dt_Q$ such that $\om:= \sum_{i=0}^{p-1} \om_i y^i$ lies in $(\pi_* \ker V_Y^n)_Q \otimes k((t_Q))$. 
\end{lemma}

\begin{proof}
As in~\cite[Lemma 5.10]{BoCaASc}, $\coef_m(h)$ denotes the coefficient of $t_Q^m$ in $h \in k((t_Q))$ and $\coef_m(\xi)$ denotes the coefficient of $t_Q^m \frac{dt_Q}{t_Q}$ in $\xi \in k((t_Q)) dt_Q$. We define the coefficients of $\om_{p-1}$ inductively to make it satisfy all the conditions above. Our goal is to construct a differential $\om_{p-1}$ that is regular at $Q$ such that $$ V_Y^n(\om)=V_Y^n\left( \sum_{i=0}^{p-1} \om_i y^i \right) = 0 . $$
Appealing to Lemma~\ref{lemVYinVX}, we obtain
$$\sum_{j\leq j_{1}\leq \ldots \leq j_n \leq p-1}
W_{j_{1},j}( W_{j_{2},j_{1}}( \ldots W_{j_{n},j_{n-1}}(\om_{j_n}) \ldots )) = 0$$
for every $0\leq j \leq p-1$. For $j=p-1$, this immediately gives $V_X^n(\om_{p-1})=0$. In the case $j<p-1$, isolating the contribution from $\om_{p-1}$ yields
\begin{equation} \label{eq5.10}
\begin{split}
\sum_{j\leq j_{1}\leq \ldots \leq p-1}
W_{j_{1},j}( W_{j_{2},j_{1}}( \ldots W_{p-1,j_{n-1}}(\om_{p-1}) \ldots )) = \\
-\sum_{j\leq j_{1}\leq \ldots \leq j_n < p-1}
W_{j_{1},j}( W_{j_{2},j_{1}}( \ldots W_{j_{n},j_{n-1}}(\om_{j_n}) \ldots )).
\end{split}
\end{equation}

First note the order of vanishing is minimized if $j_1=p-1$, in which case the powers of $f$, which has a pole at $Q$, are not diminished by the Cartier operator. The assumption on $\ord_Q(\om_i)$ then implies that on the right-hand side we have 
\begin{equation} \label{eq5.10range}
    \coef_m\left((-\sum_{j\leq j_{1}\leq \ldots \leq j_n < p-1} 
W_{j_{1},j}( W_{j_{2},j_{1}}( \ldots W_{j_{n},j_{n-1}}(\om_{j_n}) \ldots ))\right)=0
\end{equation}
for $m<-\frac{(p-1-j)d_Q}{p}$.
We write $\om_{p-1}=\sum_{m=1}^\infty b_m t_Q^m \frac{dt_Q}{t_Q}$ and put restrictions on the coefficients $b_m$. For $\om_{p-1}$ to be regular, we set $b_m=0$ for $m \leq 0$. Note that we must have $b_{m}=0$ whenever $p^n|m$, since $V_X^n(\om_{p-1})=0$. Now, we also let $b_m=0$ if $m \not \equiv 0 \bmod p^{n-1}$. We then obtain $$\om_{p-1}=\sum_{l\geq 1, p \not \; | l} b_{lp^{n-1}} t_Q^{lp^{n-1}} \frac{dt_Q}{t_Q}.$$ We now define the coefficients $b_{lp^{n-1}}$ inductively to satisfy \eqref{eq5.10}, starting with the first unspecified coefficient $b_{p^{n-1}}$. We let $0 \leq j \leq p-1$ be such that $d_Q(p-1-j) \equiv 1 \bmod p$. On the left-hand side of \eqref{eq5.10}, we have the term 
\begin{align*}
    &V_X \left( \binom{p-1}{j} \coef_{-d_Q(p-1-j)} \left((-f)^{p-1-j} \right)V_X^{n-1}\left(b_{p^{n-1}} t_Q^{p^{n-1}} \frac{dt_Q}{t_Q} \right) \right) \\
    = &V_X\left( \binom{p-1}{j} b_{p^{n-1}}^{1/p^{n-1}} \coef_{-d_Q(p-1-j)}\left( (-f)^{p-1-j} \right) dt_Q \right) \\
    = &\binom{p-1}{j} b_{p^{n-1}}^{1/p^{n}} \left(\coef_{-d_Q}\left(-f \right)\right)^{(p-1-j)/p} t^{-\frac{d_Q(p-1-j)-1}{p}-1} dt_Q.
\end{align*}
This term corresponds to setting $j_1=p-1$. This is the only term of concern, since letting $j_1<p-1$ or $l>1$ results in a greater order of vanishing at $Q$ (i.e. a smaller pole order). This implies that on the left-hand side we have  
$$\binom{p-1}{j} b_{p^{n-1}}^{1/p^{n-1}} \left( \coef_{-d_Q}\left(-f \right)\right)^{(p-1-j)/p} $$ $$ = \coef_{-\frac{d_Q(p-1-j)-1}{p}}\left(\sum_{j\leq j_{1}\leq \ldots \leq p-1}
W_{j_{1},j}( W_{j_{2},j_{1}}( \ldots W_{p-1,j_{n-1}}(\om_{p-1}) \ldots )) \right) .$$ 
Since the same coefficient of the right-hand side of \eqref{eq5.10} is determined by $\om_0, \ldots, \om_{p-2}$, which are given, we can solve for the coefficient $b_{p^{n-1}}$. 

We now continue inductively. Let $N$ be such that $b_{lp^{n-1}}$ has been defined for $l<N$. Let $j$ be such that $d_Q(p-1-j) \equiv N \bmod p$. Then the $\frac{N-d_Q(p-1-j)}{p}$-th coefficient of the left-hand side of \eqref{eq5.10} is determined by $b_{Np^{n-1}}$ and coefficients that have already been defined. Thus we can solve for $b_{Np^{n-1}}$ to arrange that the $\frac{N-d_Q(p-1-j)}{p}$-th coefficient is the same on both sides of \eqref{eq5.10}. In this way Equation \eqref{eq5.10} will be satisfied for every $j$; \eqref{eq5.10range} deals with the case $m<-\frac{(p-1-j)d_Q}{p}$. For $m>-\frac{(p-1-j)d_Q}{p}$, we let the coefficients agree by solving for the coefficient $b_{(pm+d_Q(p-1-j))p^{n-1}}$.

\end{proof}

\begin{proposition} \label{propbrpsurj}
For any $0 \leq j \leq p-1$ and $Q \in S$, the map $g'_{j+1,Q}: \G_{j+1,Q} \to M_{j,Q}$ is surjective.
\end{proposition}
\begin{proof}
As in~\cite[Proposition 5.9]{BoCaASc}, it suffices to check surjectivity after tensoring with $\cO_{X,Q}^\wedge$. As $g'_{j+1,Q}$ is an $\cO_{X,Q}$-module homomorphism, it suffices to construct $(\nu_0, \ldots , \nu_{p-1})$ such that $$g'_{j+1,Q}(\nu_0,\ldots,\nu_{p-1})= 1.$$ To this end, we pick arbitrary $\om_0, \ldots,\om_{j-1},\om_{j+1}, \ldots, \om_{p-2}$ such that $\ord_Q(\om_i)\geq -n_{Q,i}$ and set $$\om_j:=s_j\left(d_Q^{-m_{Q,j}/p^n} \frac{dt_Q}{t_Q} \right).$$ 

It is a consequence of Lemma~\ref{lemsecs} that $$\ord_Q(\om_j) \geq -m_{Q,j}-1 \geq -p^{n-1}((p-1-j)d_Q-1)-1.$$
Furthermore, by the inequality $-n_{Q,i}\geq -p^{n-1}((p-1-j)d_Q-1)-1$, the conditions of Lemma~\ref{lem5.10} are met. The lemma then tells us that there exists $\om_{p-1} \in k[\![t_Q]\!]dt_Q$ such that $\om:= \sum_{i=0}^{p-1} \om_i y^i$ lies in $(\pi_* \ker V_Y^n)_Q \otimes_{\cO_{X,Q}} k((t_Q))$. We then put $(\nu_0, \ldots, \nu_{p-1}) := \phi_\eta (\om)$ and note that $\nu_i=r_i(\om_i)$ for $i\neq j$, so $\ord_Q(\nu_i) \geq -n_{Q,i}$. Since moreover $\nu_j=0$ as $r_j \circ s_j =0$ we conclude $(\nu_0, \ldots , \nu_{p-1}) \in \G_{j+1,Q}^\wedge$. Finally, by Definition~\ref{defg'brp} we check that indeed
$$ g'_{j+1,Q}(\nu_0, \ldots , \nu_{p-1})=1 \in M_{j,Q} = M_{j,Q} \otimes_{\cO_{X,Q}} \cO_{X,Q}^\wedge,$$ which finishes the proof.
\end{proof}

\subsection{Local calculations above poles of sections} \label{secsec}

We now conduct a similar analysis for points $Q \in \supp(D_j)$. This is less technical and more similar to~\cite[Section 5C]{BoCaASc}. We fix $0 \leq j \leq p-1$ and $Q \in \supp(D_j)$ and let $t_Q$ be a uniformizer of $\cO_{X,Q}$. We let $v:=(\nu_0, \ldots, \nu_{p-1}) \in \G_{p,Q}$. We set $(\om_0, \ldots, \om_{p-1}):=\iota(\nu)$ and define $\om:= \sum_{i=0}^{p-1} \om_i y^i$. Note that this construction implies $V_Y^n(\om)=0$.

Lemma~\ref{lem5.4}.(ii) yields that $\ord(\om_j)\geq -p^n$, which gives the expansion $\om_j= \sum_{i=-p^n}^\infty b_i t_Q^i dt_Q$ in $\cO_{X,Q}^\wedge$. We compute
$$ g_{j+1,Q}(\nu) \equiv \sum_{i=-p^n}^{-1} b_i t_Q^idt_Q \bmod k[\![t_Q]\!]dt_Q.$$
This allows us to define the analogous maps $c_{j,Q}$ and $g'_{j,Q}$.

\begin{definition} \label{defg'sec}
We define a map 
\begin{align*}
    c_{j,Q}: \im(g_{j+1})_Q &\to M_{j,Q} \\
    \sum_{i=-p^n}^{-1} b_i t_Q^idt_Q \bmod k[\![t_Q]\!]dt_Q &\mapsto \sum_{i=0}^{p^n-2} b_{i-p} t_Q^{p^ni} \bmod t_Q^{p^n(p^n-1)} k[\![t_Q]\!].
\end{align*}
Then $g'_{j+1,Q}: \G_{j+1} \to M_{j,Q}$ is defined to be the composite map $c_{j,Q} \circ g_{j+1,Q}$.
\end{definition}

The exponent is multiplied by $p^n$ to make sure we land in $M_{j,Q}$ from Definition~\ref{defMjQ}. As in~\ref{secbrp}, we show that $g'_{j+1,Q}$ is surjective.

\begin{lemma} \label{lem0sec}
Let $\nu \in \G_{j+1,Q}$ and $g'_{j+1,Q}(\nu)=0$. Then $\ord_Q(\om_j)\geq -n_{Q,j}=0$, meaning that $\nu \in \G_{j,Q}$.
\end{lemma}

\begin{proof}
This is a straightforward amendment of~\cite[Lemma 5.12]{BoCaASc}. First off, we have $\om_j=\nu_j+s_j(V_X^n(\om_j))$. Using Corollary~\ref{corVYker}, we expand
$$ V_X^n (\om_j) = - \sum_{i=j+1}^{p-1} \left( \sum_{j \leq j_1 \leq \ldots \leq j_{n-1} \leq i} W_{j_2,j_1}(\ldots W_{i,j_{n-1}}(\om_i))\right) = (a_1+a_2 t_Q + \ldots)dt_Q . $$
Note that regularity follows from the fact that the divisors $D_j$ are chosen to have disjoint support, so that $\ord_Q(w_i)\geq-p^n$. Then we obtain
\begin{equation} \label{eqomjsec}
\om_j=\nu_j+a_1^{p^n}s_j(dt_Q) + a_2^{p^n}dt_Q^{p^n}s_j(dt_Q) + \ldots.
\end{equation}
By the construction of $s_j$, we have $\ord_Q(s_j(dt_Q)) \geq -p^n$, whence $\ord_Q(\om_j) \geq -p^n$. Moreover, since $g'_{j+1,Q}(\nu)=0$, the coefficient of $t_Q^i$ in $\om_j$ is zero for $-p^n \leq i \leq -2$. Finally, the coefficient of $t_Q^{-1}$ is zero in $\nu_j$, and therefore in $\om_j$, since $V_X^n(\nu_j)=0$. We conclude that $\ord_Q(\om_j) \geq 0$, as required.
\end{proof}

\begin{proposition} \label{propDjsurj}
For $0 \leq j \leq p-1$ and $Q \in \supp(D_j)$, the map $g'_{j+1,Q}: \G_{j+1,Q} \to M_{j,Q}$ is surjective.
\end{proposition}

\begin{proof}
We show that an arbitrary element $\sum_{i=0}^{p^n-2} c_i t_Q^{p^n i} \in M_{j,Q}$ lies in the image of $g'_{j+1,Q}$. We simply put $\nu_j:= \sum_{i=0}^{p^n-2} c_i t_Q^{i-p^n}$ and $\nu_i=0$ for $i \neq j$. Then, as $\nu_j$ has no $t_Q^{-1}$ term, $V_X^n(\nu_j)=0$ and hence $\nu \in \G_{j+1,Q}$. By construction we have $g'_{j+1,Q}(\nu)=\sum_{i=0}^{p^n-2} c_i t_Q^{p^n i}$, as desired.
\end{proof}

\subsection{Local calculations above \texorpdfstring{$Q'$}{}}

Finally we account for the behavior at the auxiliary point $Q'$. As we will now see, this amounts to computing the quotient sheaf $\G_0/\G_{-1}$, where $\G_{-1}:=\im(\phi)=\pi_* \ker(V_Y)$. As before, $t_{Q'}$ denotes a uniformizer of $\cO_{X,Q'}$. We define $S_{-1}:=\{Q'\}$.

\begin{lemma} \label{lemMjQ'}
For $0 \leq j \leq p-1$, the inclusion $\G_{j,Q'} \to \G_{j+1,Q'}$ is an isomorphism.
\end{lemma}
\begin{proof}
We show surjectivity. Let $\nu \in \G_{p,Q}$ and put $\iota(\nu)=: (\om_0, \ldots, \om_{p-1})$. By Lemma~\ref{lem5.4}.(iii) we have $\ord_{Q'}(\om_i) \geq -p^{n-1}d_{Q'}(p-1-i)$ for all $i$, meaning that $\nu \in G_{0,Q}$. This implies in particular that $\G_{j,Q'} \to \G_{j+1,Q'}$ is surjective. 
\end{proof}

On the other hand, the sheaves $\G_{-1}$ and $\G_0$ are only different at the point $Q'$. We will now construct the sheaf $M_{-1}$ and show it is isomorphic to the quotient sheaf $\G_0/\G_{-1}$.

\begin{definition} \label{defM-1}
Define $M_{-1,Q'}:=\bigoplus_{i=0}^{p-1} \left(k[\![t_{Q'}^{p^n} ]\!]/(t_{Q'}^{p^{n-1}d_{Q'}(p-1-i)})  \right)^{\oplus (p^n-1)}$.
We let \\ $M_{-1}:=Q_*(M_{-1,Q'})$, the skyscraper sheaf supported at $Q'$ with stalk $M_{-1,Q'}$. 
\end{definition}

Note that $M_{-1,Q'}$ is well-defined since $d_{Q'}$ is divisible by $p$.

\begin{proposition} \label{propQ'coker}
The cokernel of the natural inclusion $\G_{-1} \to \G_0$ is isomorphic to $M_{-1}$.
\end{proposition}

\begin{proof}
The strategy of~\cite[Lemma 5.16]{BoCaASc} can largely be copied. This approached is based on switching to a model $(y')^p-y'=f'$ such that $f'$ does not have a pole at $Q'$. One then studies
\begin{align*}
N_{-1}&:=\left\{ \om=\sum_{i=0}^{p-1} \om'_i (y')^i : \om \in (\pi_* \ker V_Y^n)_{Q'} \text{ and } \om'_i \in \Om^1_{X,Q'}   \right\} \cong \G_{-1,Q'} \\
N_0 &:= \left\{ \om=\sum_{i=0}^{p-1} \om'_i (y')^i : \om \in (\pi_* \ker V_Y^n)_{Q'} \text{ and } \om'_i \in t_{Q'}^{-p^{n-1}d_{Q'}(p-1-i)} \Om^1_{X,Q'}   \right\} \cong \G_{0,Q'} 
\end{align*}
Analyzing the completions of these modules, using a modified version of Proposition~\ref{propphieta}, shows
\begin{align*}
    N_{-1}^\wedge &\cong \bigoplus_{i=0}^{p-1} (\pi_* \ker V_X^n)_{Q'}^\wedge \\
    N_0^\wedge &\cong \bigoplus_{i=0}^{p-1} t_{Q'}^{-p^{n-1}d_{Q'}(p-1-i)} (\pi_* \ker V_X^n)_{Q'}^\wedge.
\end{align*}
As a $k[\![t_Q^{p^n}]\!]$-module, $(\ker V_X^n)_{Q'}^\wedge$ is generated by $dt_{Q'}, \ldots , t_{Q'}^{p^n-2} dt_{Q'}$. We exclude precisely the exponents that are $-1 \bmod p^n$ as these will survive $V_X^n$. We then write
\begin{align*}
    \G_{-1}/\G_0 &= N_{-1}/N_0 = \frac{\bigoplus_{i=0}^{p-1} (\pi_* \ker V_X^n)_{Q'}^\wedge}{\bigoplus_{i=0}^{p-1} t_{Q'}^{-p^{n-1}d_{Q'}(p-1-i)} (\pi_* \ker V_X^n)_{Q'}^\wedge}  \\
    &\cong \bigoplus_{i=0}^{p-1} (\ker V_X^n)_{Q'}^\wedge/t_{Q'}^{p^{n-1}d_{Q'}(p-1-i)} (\ker V_X^n)_{Q'}^\wedge = M_{-1,Q'},
\end{align*}
which finishes the proof.
\end{proof}

\subsection{Short exact sequences}

In the preceding sections, we have computed the stalks of $\G_{j+1}/\G_j$ at ramified points, poles of sections and the auxiliary point $Q'$. Using these computations, we can now completely describe the associated graded sheaves of the Galois invariant filtration presented in \eqref{eqGfilt}. We put the maps of stalks together to form a maps of sheaves.

\begin{definition} \label{defg'global}
Let $0 \leq j \leq p-1$. Recall $S_j = S \cup \supp(D_j)$. Define the map
$$
\begin{tikzcd}
c_j: \im(g_{j+1}) \arrow[r] & \bigoplus_{Q \in S_j} Q*(\im(g_{j+1})_Q) \arrow[r, "{\bigoplus_{Q \in S_j} Q_*(c_{j,Q})}"] & {\bigoplus_{Q \in S_j} Q_*(M_{j,Q}) = M_j.}
\end{tikzcd}
$$
Then define $g'_{j+1}$ to be the composition $c_j \circ (g_{j+1}|_{\G_{j+1}})$. Moreover, define $g'_0$ to be the natural projection $\G_0 \to \G_0/\G_{-1}$ composed with the isomorphism from Proposition~\ref{propQ'coker}.
\end{definition}

Note that this definition of the global map $g'_{j+1}$ indeed localizes to the local map defined in Definitions~\ref{defg'brp} and~\ref{defg'sec}.
This definition allows us to formulate the exact sequence relating $\G_j$ to $\G_{j+1}$.

\begin{theorem} \label{thmexseqMj}
For $-1 \leq j \leq p-1$, there is an exact sequence
\begin{equation} \label{exseqMj}
\begin{tikzcd}
0 \arrow[r] & \G_j \arrow[r] & \G_{j+1} \arrow[r, "g'_{j+1}"] & M_j \arrow[r] & 0 .
\end{tikzcd}
\end{equation}
\end{theorem}

\begin{proof}
For $j=-1$, this is the content of Proposition~\ref{propQ'coker}, so we assume $j \geq 0$. 

Note that exactness on the left is a direct consequence of the construction of the sheaves $\G_j$. We check exactness in the middle and on the right locally. 

First we localize at the points $Q \notin S_j$. There $M_{j,Q}=0$ by definition and $\G_{j,Q} \to \G_{j+1,Q}$ is an isomorphism by Lemma~\ref{lemGjbasics}.(ii), yielding exactness at $Q$. 

We now localize at a branch point $Q \in S$. Proposition~\ref{propbrpsurj} provides exactness on the right. The observation that $c_{j,Q}$ is injective then implies that $\ker(g'_{j+1})=\ker(g_{j+1})$, which proves exactness in the middle when combined with Lemma~\ref{lemGjbasics}.(i).

Finally, take a point $Q \in \supp(D_j)$. Exactness on the right is the content of Proposition~\ref{propDjsurj}. Exactness in the middle follows from Lemma~\ref{lemGjbasics}.(ii) and Lemma~\ref{lem0sec}.
\end{proof}

The remainder of the paper will revolve around bounding the dimensions $H^0(X,M_j)$ and thereby bounding the $n$-th $a$-number $a_Y^n=H^0(X,\G_{-1})$.

\begin{example} \label{eg7again}
We revisit Example~\ref{eg7} using the language acquired in this section. Recall that the curve $Y_2$ is given by $y^7-y=t^{-4}+t^{-3}$. As we saw, the global section of interest was
$$ \nu_{4,2}:= (0, 0, 0, 0, t^{-2}dt, 0, 0) \in \bigoplus_{i=0}^{p-1} H^0(\PP^1,\ker V_{\PP^1}^2(F_*^2 E_i)) = H^0(\PP^1, \G_7).$$
Since $X=\PP^1$, only the analysis in section~\ref{secbrp} is needed, at the ramified point $Q=0$ with $t_Q=t$. By the definition of $\ph$, $\nu_{4,2}$ lies in $H^0(\PP^1,\G_4)$. In Example~\ref{eg7} we computed $$\om := \ph_\eta^{-1}(\nu_{4,2})=t^{-2}dt y^4 + t^{-50}dt.$$ Hence $\nu_{4,2}$ lies in $H^0(\PP^1,\G_1)$, but not in $H^0(\PP^1,\G_0)$, because $\om_0=t^{-50}dt$ is not a global section of $F_*^n(\Om_X^1(E_i))$ (see Definition~\ref{defGj}). As $g'_0$ extracts the irregular terms of $\om_0$. we have $H^0(\PP^1,M_{0,Q})=\spn_k\{t^{-50}\}$. Finally, we verify that $H^0(\PP^1,\G_{-1})=H^0(\PP^1,\G_0)$ and $H^0(\PP^1,\G_1)=H^0(\PP^1,\G_7)$, so we have $a^2_{Y_2}=9-\dim_k H^0(\PP^1,M_{0,Q})=8$.

\end{example}

\section{Bounds} \label{sec6}

As before, let $\pi:Y\to X$ be an Artin-Schreier cover of curves. In this section we use Theorem~\ref{thmexseqMj} to bound $$a_Y^n=\dim_k H^0(Y,\ker V_Y^n) = \dim_k H^0(X,\G_{-1})$$ in terms of $a_X^n$ and the ramification invariants of $\pi$.

\subsection{Abstract bounds}
In this subsection we formulate some abstract bounds, which we express in more concrete quantities later in the section. As in~\cite[Definition 6.1]{BoCaASc}, we first encode the information of the exact sequences \eqref{exseqMj} into one linear map.

\begin{definition} \label{deftgU}
For $0 \leq j \leq p-1$, let $\Tilde{g}_{j+1}: H^0(X,\G_p) \to H^0(X,M_j)$ be the map on global sections induced by $c_j \circ g_{j+1}$, as defined in Definitions~\ref{defg'brp},~\ref{defg'sec} and~\ref{defggj}. Define
\begin{align*}
    \Tilde{g}:H^0(X,\G_p) &\to \bigoplus_{j=0}^{p-1} H^0(X,M_j) \\
    \nu &\mapsto (\Tilde{g}_1(\nu),\ldots ,\Tilde{g}_p(\nu)).
\end{align*}
Let $\Tilde{g}_0: H^0(X,\G_0) \to H^0(X,M_{-1})$ be the map of global sections induced by $g'_0$ from Definition~\ref{defg'global}.
Finally, set 
\begin{align*}
N_1(X,\pi)&:=\dim_k \im (\tg) \\ 
N_2(X,\pi)&:=\dim_k \im(\tg_0). \\
N(X,\pi)&:= N_1(X,\pi) + N_2(X,\pi) \\
U^n(X,\pi)&:= \sum_{i=0}^{p-1} \dim_k H^0(X,\ker V_X^n(F_*^n(E_i+p^nD_i))) - N(X,\pi),
\end{align*}
\end{definition}

\begin{proposition} \label{propU(X,pi)}
We have $a_Y^n = U^n(X,\pi)$.
\end{proposition}

\begin{proof}
This is not fundamentally different from Lemma 6.3 and Lemma 6.5 in~\cite{BoCaASc}. 

The statement $\ker(\tg) = H^0(X, \G_0)$ can be made intuitive as follows. If $\nu\in H^0(X,\G_p)$ fails to be in $H^0(X, \G_0)$, then $\om:=\iota(\nu)$ fails to be regular at a point $Q \neq Q'$. In that case this failure is recorded by $M_j$, by construction. 

It can be proved rigorously by applying Theorem~\ref{thmexseqMj} repeatedly. For $\nu \in H^0(X,\G_0)$, the left exactness of $H^0$ applied to \eqref{exseqMj} with $j=0$ implies that $\tg_1'(\nu)=0$. By increasing induction on $j$, we find that $\tg_j'(\nu)=0$ for every $1 \leq j \leq p$ and hence $\nu \in \ker (\tg)$. For the reverse inclusion suppose $\nu \in \ker (\tg)$, so that $\nu \in H^0(X,\G_p)$ and $\tg_j'(\nu)=0$ for each $j$. Now, taking $j=p-1$ in \eqref{exseqMj} implies that $\nu \in H^0(X,\G_{p-1})$. Continuing the descending induction on $j$ yields that $\nu \in H^0(X,\G_0)$, as desired. Then the rank-nullity theorem yields
\begin{equation} \label{eqN1}
\begin{aligned} 
    \dim_k H^0(X,\G_0)&=\dim_k H^0(X,\G_p)-\dim_k \im (\tg) \\
    &= \sum_{i=0}^{p-1} \dim_k H^0(X,\ker V_X^n(F_*^n(E_i+p^nD_i))) - N_1(X,\pi).
\end{aligned}
\end{equation}

Finally, Proposition~\ref{propQ'coker} and the left exactness of $H^0$ provide that $$\ker( \tg_0)=H^0(X,\G_{-1})=H^0(X,\pi_*\ker V_Y^n).$$ Again applying the rank-nullity theorem and \eqref{eqN1} yields
\begin{align*}
    a_Y^n &= \dim_k H^0(X,\G_{-1}) \\
    &= \dim_k H^0(X,\G_0) - \dim_k \im(\tg_0) \\
    &= \sum_{i=0}^{p-1} \dim_k H^0(X,\ker V_X^n(F_*^n(E_i+p^nD_i))) - N_1(X,\pi) - N_2(X,\pi) \\
    &= U^n(X,\pi).
\end{align*}
\end{proof}

To obtain an upper bound on $a_Y^n$, we will bound $U^n(X,\pi)$ from above. We also introduce an abstract lower bound $L^n(X,\pi)$. Recall that we have fixed a generator $\tau$ of $\Gal(k(Y)/k(X)) \cong \Z/p\Z$ that acts on $\pi_* \Om_Y^1$ and preserves $\G_j$ for each $j$.

\begin{definition} \label{defL(X,pi)}
Let $i$ and $j$ be integers between $0$ and $p-1$. We make the following definitions:
\begin{itemize}
    \item $U_j:=\ker\left( (\t-1)^{j+1}: H^0(X,\G_{j+1}) \to H^0(X,\G_{j+1})  \right).$
    \item For $Q \in S$, $c(i,j,Q)$ is the number of integers $m$ that are $-1 \bmod p^n$ such that $-n_{Q,j}-d_Qp^{n-1}(j-i) \leq m < -n_{Q,i}$.
    \item For $Q \in \supp(D_j)$, $c(i,j,Q):=p^n-1$.
    \item $c(i,j,Q'):=\frac{d_{Q'}(p^n-1)(p-l-i)}{p}$.
    \item $L^n(X,\pi):= \max_{0 \leq j \leq p-1} \left( \dim_k U_j - \sum_{i=0}^j \sum_{Q \in S_i'} c(i,j,Q)  \right)$.
\end{itemize}
\end{definition}

The reason for this definition is that we will be able to compute $\dim_k U_j$, while the $c(i,j,Q)$'s together give an upper bound on $\dim_k \tg (U_j)$, such that the rank-nullity theorem yields the lower bound $L^n(X,\pi)$ for $a_Y^n$. 

\begin{lemma} \label{lemUj}
There is an isomorphism $$U_j \cong \bigoplus_{i=0}^j H^0(X, \ker V_X^n(F_*^n(E_i+p^nD_i))).$$
\end{lemma}

\begin{proof}
As in~\cite[Lemma 6.7]{BoCaASc}, both are isomorphic to 
$$ \left \{ \om=\sum_{i=0}^j\om_i y^i \in (\pi_* \ker V_Y)_\eta \, : \, \ph_\eta(\om) \in H^0(X,\G_p)  \right \},$$ since $\tau-1$ reduces the maximum power of $y$.
\end{proof}

We now show how the $c(i,j,Q)$ can be used to bound the image of $\tg$ from above. For $-1 \leq i \leq j \leq p-1$, define 
$$\psi^j_{i+1}:= H^0(g_{i+1}')|_{U_j \cap H^0(X,\G_{i+1})}: U_j \cap H^0(X,\G_{i+1}) \to H^0(X,M_i) = \bigoplus_{Q \in S_i'} M_{i,Q}.$$
Let $\psi^j_{i+1,Q}$ be $\psi^j_{i+1}$ composed with the projection to $M_{i,Q}$.

In Proposition~\ref{propL(X,pi)} it will be crucial to bound $\dim_k \im(\psi^j_{i+1,Q})$ from above, which the following lemma does.

\begin{lemma} \label{lemcijQ}
We have $\dim_k \im(\psi^j_0)=\sum_{i=0}^j c(i,j,Q')$. For $0 \leq i \leq j \leq p-1$ and $Q \in S_i$, we have $\dim_k \im(\psi^j_{i+1,Q}) \leq c(i,j,Q)$. 
\end{lemma}

\begin{proof}
For the first statement, recall that $M_{-1}$ is only supported at $Q'$. By the proof of Proposition~\ref{propQ'coker}, the image of $H^0(g_0)$ when restricted to $$U_j \cap H^0(X,\G_{j+1})= H^0(X,\ker((\tau-1)^{j+1}: \G_0 \to \G_0))$$ has dimension at most
$$\sum_{i=0}^j \frac{p^{n-1}d_{Q'}(p^n-1)(p-1-i)}{p^n}= \sum_{i=0}^j c(i,j,Q'),$$
as desired.

For the second statement, consider $Q \in \supp(D_j)$ first. By Definition~\ref{defg'sec}, $\im(\psi^j_{i+1,Q})$ has dimension at most $p^n-1$, which we defined $c(i,j,Q)$ to be. 

On the other hand, suppose $Q \in S$. For $\nu \in U_j$, write 
\begin{align*}
    \nu&=(\nu_0, \ldots, \nu_j, 0, \ldots, 0), \\
    \intertext{such that} 
    \iota(\nu)&=:(\om_0, \ldots , \om_j, 0 \ldots, 0) \\
    \intertext{where}
    \om_i&=\nu_i + s_i\left( - \sum_{l=i+1}^{j} \left( \sum_{i \leq j_1 \leq \ldots \leq j_{n-1} \leq l} W_{j_1,j}(\ldots W_{l,j_{n-1}}(\om_l))\right)  \right)
\end{align*}
by \eqref{eqphiinv}. Lemma~\ref{lem5.4} then provides
\begin{equation} \label{eqomiordQ}
\ord_Q(\om_i) \geq -n_{Q,j} - p^{n-1}d_Q(j-i).
\end{equation}

By Definition~\ref{defg'brp}, $\psi^j_{i+1,Q}(\nu)$ records the coefficients in the local expansion of $\om$ at $Q$ whose exponent is $-1 \bmod p^n$ and smaller than $-n_{Q,i}$. By \eqref{eqomiordQ}, this exponent is at least $-n_{Q,j}-p^{n-1}d_Q(j-i)$. This implies that the dimension of $\im(\psi^j_{i+1,Q})$ is at most $c(i,j,Q)$.

\end{proof}

\begin{proposition} \label{propL(X,pi)}
We have $a_Y^n \geq L^n(X,\pi)$.
\end{proposition}

\begin{proof}
For fixed $j$, this is proved using descending induction on $i$, going from $\G_{j+1}$ to $\G_{-1}\cong \pi_* \ker V_Y^n$. Using Lemma~\ref{lemUj}, we observe $U_j \subseteq H^0(uX, \G_{j+1})$ and hence
\begin{equation*} \label{eqdimUjGistart}
\dim_k (U_j \cap H^0(X,\G_{j+1})) = \dim_k U_j.
\end{equation*}
Then, by Theorem~\ref{thmexseqMj} and the left exactness of $H^0$ we obtain for each $0 \leq i \leq j$ the exact sequence
\begin{equation} \label{exseqUjpsi}
    0 \to U_j \cap H^0(X,\G_i) \to U_j \cap H^0(X,\G_{i+1}) \to \im(\psi^j_{i+1}) \to 0.
\end{equation}
Furthermore, the inclusion 
$$ \im(\psi^j_{i+1}) \subseteq \bigoplus_{Q\in S_i} \im(\psi^j_{i+1,Q}),$$
together with Lemma~\ref{lemcijQ} gives the bound 
$$ \dim_k \im(\psi^j_{i+1}) \leq \sum_{i=0}^j c(i,j,Q).$$
We now apply the rank-nullity theorem to \eqref{exseqUjpsi} to obtain
$$ \dim_k H^0(X,\G_i) \geq H^0(X,\G_{i+1}) - \sum_{Q \in S_i} c(i,j,Q).$$
Iterating this inequality for descending $i$ yields
$$ \dim_k H^0(X,\G_0) \geq \dim_k U_j - \sum_{i=0}^j \sum_{Q\in S_i} c(i,j,Q).$$
Finally, applying the rank-nullity theorem again for $i=-1$ means the involvement of the auxiliary point $Q'$. Appealing again to Lemma~\ref{lemcijQ} gives the desired inequality
$$a_Y^n = \dim_k H^0(X,\G_{-1}) \geq \dim_k U_j - \sum_{i=0}^j \sum_{Q\in S_i'} c(i,j,Q). $$
\end{proof}

\subsection{Tools}

In this subsection, we make the abstract bounds more explicit by computing the dimension of certain spaces of global sections. The main ingredient for these computations is Tango's theorem~\cite[Theorem 15]{Tango}.

\begin{proposition} \label{proptango}
For a curve $X$ over $k$, define the \emph{Tango number} by
$$m(X):= \max \left \{ \sum_{P\in X(k)} \left\lceil \frac{\ord_P(df)}{p} \right\rceil \; | \; f \in k(X) - k(X)^p \right\}.$$
For any line bundle $\mathcal{L}$ with $\deg \mathcal{L} > m(X)$, the map 
$$ V_X^n: H^0(X,\Om_X^1 \otimes \mathcal{L}^p) \to H^0(X,\Om_X^1 \otimes \mathcal{L}) $$
is surjective. 
\end{proposition}

\begin{proof}
We prove by induction on $n$. For $n=1$, this is the statement of~\cite[Theorem 15]{Tango}. For $n>1$, we apply Tango's theorem to the line bundle $\mathcal{L}^{p^{n-1}}$, which has degree $p^{n-1}\deg \mathcal{L} \geq \deg \mathcal{L}>m(X)$, since $m(X)\geq -1$. Thus by induction the composition
$$
\begin{tikzcd}
{V_X^n:H^0(X,\Omega_X^1 \otimes \mathcal{L}^{p^n})} \arrow[r, "V_X"] & {H^0(X,\Omega_X^1 \otimes \mathcal{L}^{p^{n-1}})} \arrow[r, "V_X^{n-1}"] & {H^0(X,\Omega_X^1 \otimes \mathcal{L})}
\end{tikzcd}
$$
is surjective.
\end{proof}

For the next proof, we introduce the abbreviation 
$$\d(H^0(X,\Om_X^1(E))) := \dim_k \ker \left( V_X^n: H^0(X,\Om_X^1(E)) \to H^0(X,\Om_X^1(E)) \right).$$
Recall that $a_X^n=\d(H^0(X,\Om_X^1))$.

\begin{corollary} \label{corwrong}
Let $D$ and $R$ be divisors on $X$, with $R=\sum_i r_i P_i$ and $0\leq r_i <p^n$. Under the condition $\deg(D) > \max(m(X),0)$, we have
$$\d (H^0(X,\Om_X^1(p^nD+R))) = (p^n-1) \deg (D) + \sum_i \left( r_i - \left\lceil \frac{r_i}{p^n} \right\rceil \right). $$
For general $D$ we have the inequality
$$0 \leq \d(H^0(X,\Om_X^1(p^nD+R))) - \left( (p^n-1)\deg(D) + \sum_i \left(r_i - \left\lceil \frac{r_i}{p}\right\rceil \right)\right) \leq a_X^n.$$
\end{corollary}

\begin{proof}
For the first assertion, we start with the case $R=0$. Applying Proposition~\ref{proptango} to the case $\mathcal{L}=\cO_X(D)$ yields surjectivity of $V_X^n$. Then the Riemann-Roch theorem implies
\begin{align*}
    \d(H^0(X,\Om_X^1(p^nD))) &= \dim_k H^0(X,\Om_X^1(p^nD)-\dim_k H^0(X,\Om_X^1(D)) \\
    &= (g_X-1)+p^n \deg(D) - \left( (g_X-1) + \deg(D) \right) \\
    &= (p^n-1) \deg(D)
\end{align*}
as desired. We now suppose $R \neq 0$ and prove the general equality by induction on the cardinality of $\supp(R)$. As in~\cite[Corollary 6.13]{BoCaASc}, we have for any closed point $P \in X \setminus \supp(R)$ and for any divisor $E$ on $X$ with $\deg(E)\geq 0$ the inequality
$$ 0 \leq \d(H^0(X,\Om_X^1(p^nD+R))) - \d(H^0(X,\Om_X^1(p^nD+R+[P]))) \leq 1.$$
where this difference is zero if $p^n | \ord_P(p^nD+R)$. Since we have
\begin{align*}
\d(H^0(X,\Om_X^1(p^nD+p^n[P]+R)&=\d(H^0(X,\Om_X^1(p^nD+R)))+(p^n-1)
\intertext{by our induction step and the assumption $\deg(D)>\max(m(X),0)$, it follows that}
\d(H^0(X,\Om_X^1(p^nD+R+rP)))&=\d(H^0(X,\Om_X^1(p^nD+R)))+(r-1),
\end{align*}
for every $0<r<p^n$, which finishes the proof of the equality.

We now prove the inequality. For the upper bound, we pass from $a_X^n=\d(H^0(X,\Om_X^1))$ to $\d(H^0(X,\Om_X^1(p^nD+R)$, such that the number of times the dimension can increase by $1$ is at most $(p^n-1)\deg(D)+\sum (r_i-\lceil r_i/p^n\rceil)$. This establishes the upper bound. 

For the lower bound, note that the argument in~\cite[Corollary 6.13]{BoCaASc} does not work. Instead, the following argument, provided by Jeremy Booher, does work and extends to the case $n>1$. We fix a divisor $E$ on $X$ with $\deg(E)>\max(m(X),0)$ and $p^nE \geq p^nD+R$. We pass from $$\d(H^0(X,\Om_X^1(p^nE)))=(p^n-1)\deg(E) $$ to $\d(H^0(X,\Om_X^1(p^nD+R))).$ When we subtract a point $P$ from the divisor, the dimension of the kernel of $V_X^n$ could decrease by $1$, except if $p^n$ divides the order at $P$ of the resulting divisor. Hence there are at most $(p^n-1)\deg(E-D)-\sum (r_i-\lceil r_i/p\rceil)$ times the dimension decreases by $1$, establishing the lower bound
$$ \d(H^0(X,\Om_X^1(p^nD+R))) \geq (p^n-1)\deg(D) + \sum (r_i -\lceil r_i/p \rceil).$$
\end{proof}

In the proof of Theorem~\ref{thmmain}, Corollary~\ref{corwrong} will be crucial for computing $\dim_k U_j$ and therefore the lower bound $L^n(X,\pi)$. Moreover, it is used to bound $N(X,\pi)$ and therefore the lower bound $U^n(X,\pi)$.

\subsection{Estimating \texorpdfstring{$N(X,\pi)$}{}}

In order to bound $U^n(X,\pi)$, as defined in Definition~\ref{deftgU}, from above, we bound $N(X,\pi)$ from below. Essentially, we investigate how the (local) leading term of $f$ contributes to $\im(\tg)$ and $\im(\tg_0)$. This explains why it is shown in Corollary~\ref{corsharp} that the maximal value of $a_Y^n$ is realized when $f=t^d$ if $X$ is the projective line. We investigate $N(X,\pi)$ using set of the triples we now define. First, fix an ordering on $\bigcup S_j'$ with $Q'$ being the smallest element.

\begin{definition} \label{defQlj}
Let $T$ be the set of triples $(Q,l,j)$ such that
\begin{itemize}
    \item $Q \in S_j'$ and $0 \leq j \leq p-1$.
    \item $0 < l \leq \ord_Q(E_j+p^nD_j)$.
    \item $l \not \equiv 1 \bmod p^n$.
    \item If $Q\in S$, then the first $n$ digits of the $p$-adic expansion of $(1-l)/d_Q$ sum up to at most $j$.
\end{itemize}
\end{definition}

\begin{definition} \label{defexists}
For a triple $(Q,l,j) \in T$, let $\xi_{Q,l,j}$ be a differential in the kernel of $V_X^n$ on $$ H^0\left(X,\Om_X^1\left(\sum_{P<Q} \ord_P(E_j+p^nD_j)P + (l-1)Q\right) \right)$$ but not in the kernel of $V_X^n$ on
$$ H^0\left(X,\Om_X^1\left(\sum_{P<Q} \ord_P(E_j+p^nD_j)P + lQ\right) \right)$$ if such a differential exists. If $\xi_{Q,l,j}$ exists, we define
$$ \nu_{Q,l,j}:= (0, \ldots ,0,\xi_{Q,l,j},0,\ldots ,0) \in H^0(X,\G_{j+1})$$ with $\xi_{Q,j,l}$ in the $j$-th component. We then say that $\nu_{Q,l,j}$ exists.
\end{definition}

Note that $\xi_{Q,l,j}$ has the following properties:
\begin{itemize}
    \item $V_X^n(\xi_{Q,l,j})=0$.
    \item $\ord_P(\xi_{Q,l,j}) \geq -\ord_P(E_j+p^nD_j)$ for $P<Q$
    \item $\ord_Q(\xi_{Q,l,j}) = -l$.
    \item $\ord_{P}(\xi_{Q,l,j}) \geq 0$ for $P>Q$.
\end{itemize}

The remainder of this subsection will be devoted to proving the following proposition.

\begin{proposition} \label{propNTB}
We have $N(X,\pi) \geq \#T - B$, where $B$ is the number of triples $(Q,l,j)$ such that $\nu_{Q,l,j}$ doesn't exist.
\end{proposition}

The proof is based on the linear independence of the elements $\tg(\nu_{Q,l,j})$ for $Q \in S_j$ and $\tg_0(\nu_{Q',l,j})$, bounding $N_1(X,\pi)$ and $N_2(X,\pi)$, respectively, from below.

\begin{definition}
We order the set $T$ as follows. We set $(Q,l,j) < (Q',l',j')$ in either of the following cases:
\begin{itemize}
    \item $j<j'$;
    \item $j=j'$ and $Q<Q'$;
    \item $j=j'$, $Q=Q'$ and $l<l'$.
\end{itemize}
\end{definition}

For a triple $(Q,l,j) \in T$, we define 
$$U_{Q,l,j}:= \spn_k \left\{ \tg(\nu_{P,a,i}) \; : \; T \ni (P,a,i) < (Q,l,j) \text{ and $\nu_{P,a,i}$ exists} \right\}. $$

The following lemmas show that each triple $(Q,l,j)$ increases $N(X,\pi)$, provided $\nu_{Q,l,j}$ exists, treating the three cases $Q\in S$, $Q\in \supp(D_j)$ and $Q=Q'$.

\begin{lemma} \label{lemQljbrp}
Let $(Q,l,j) \in T$ with $Q\in S$ and suppose $\nu_{Q,l,j}$ exists. Then $\tg(\nu_{Q,l,j}) \notin U_{Q,l,j}$.
\end{lemma}

\begin{proof}
Recall that $\tg(\nu_{Q,l,j})$ lies in $\bigoplus_{i=0}^{p-1} H^0(X,M_i)$. Let $\tg_{m+1,Q}$ denote the composition of $\tg_{m+1}$ with projection to the coordinate corresponding to $Q$. Let 
$$ (\om_0, \ldots , \om_j, 0, \ldots ,0) := \iota(\nu_{Q,l,j}),$$
such that $\ph_\eta(\sum_i \om_i y^i) = \nu_{Q,l,j}$. By Proposition~\ref{propphieta}, we have $\om_j=\xi_{Q,l,j}$ and
\begin{equation} \label{eqom_m}
    \om_m=s_m\left(  - \sum_{i=m+1}^{p-1} \left( \sum_{m \leq j_1 \leq \ldots \leq j_{n-1} \leq i} W_{j_1,m}(\ldots W_{i,j_{n-1}}(\om_i) \ldots )\right) \right)
\end{equation}
for $m<j$. 

After scalar multiplication if necessary, we write $\xi_{Q,l,j}$ locally at $Q$ as
$$\xi_{Q,l,j} = (t_Q^{-l} + \ldots ) dt_Q.$$
Since $p \not | d_Q$, we have $(1-l)/d_Q \in \Z_p$. It has the expansion
$$\frac{1-l}{d_Q} = c_0 + c_{1} p + \ldots + c_{n-1} p^{n-1} + \ldots,$$
such that
\begin{equation}
     p^n \; | \; l-1+d_Q(c_0+ c_1p + \ldots + c_{n-1} p^{n-1}).  \label{eq1-lpexp}
\end{equation} 
By Definition~\ref{defQlj}, we have $\sum_{r=0}^{n-1} c_r \leq j$. This allows us to specialize the sum in \eqref{eqom_m} to the indices 
\begin{align*}
i&=j \\
j_{n-1} &= j-c_0 \\
&\vdots \\
j_1&=j_2-c_{n-2} \\
m&=j_1-c_{n-1}.
\end{align*}
By virtue of \eqref{eq1-lpexp}, this term is then non-zero. Indeed, we obtain 
$$ l + d_Q(j-j_{n-1}) \equiv 1 \bmod p,$$ such that $W_{j,j_{n-1}}(\xi_{Q,l,j})$ is non-zero and has order $-(l+d_Q(j-j_{n-1})-1)/p-1$ at $Q$. For the next step, we obtain
$$ (l+d_Q(j-j_{n-1})-1)/p+1 + d_Q(j_{n-1}-j_{n-2}) \equiv 1 \bmod p$$
and so forth. By Lemma~\ref{lemsecs}, this assures that $\om_m$ has a non-zero term of order 
\begin{equation} \label{ordQomm}
\ord_Q(\om_m) = l-1+d_Q(c_0+ c_1p + \ldots + c_{n-1} p^{n-1})-1,
\end{equation}
which $\tg_{m+1,Q}(\nu_{Q,l,j})$ records. To finish the proof, we now show that this term is not realized by a triple $(P,a,i)<(Q,l,j)$. In that case we have $\ord_Q(\xi_{P,a,i}) > \ord_Q(\xi_{Q,l,j})=-l$ as in~\cite[Lemma 6.20]{BoCaASc}. Note that, given $m$, the solution $j_1 \leq p-1$ is unique (if it exists). Similarly, the solutions $j_2, \ldots ,j_{n-1}, j$ are all unique. Moreover, non-leading terms of $f$ or $\xi_{P,a,i}$ will always increase the order at $Q$. Thus an order at $Q$ as low as \eqref{ordQomm} will not be reached for $(P,a,i)<(Q,l,j)$. This implies that $\tg_{m+1,Q}(\nu_{Q,l,j})$ is not in the span of $\{\tg_{m+1,Q}(\nu_{P,a,i}) \; : \; (P,a,i) < (Q,l,j) \}$ and hence $\tg(\nu_{Q,l,j}) \notin U_{Q,l,j}$.
\end{proof}

\begin{lemma} \label{lemQljDj}
Let $(Q,l,j) \in T$ with $Q \in \supp (D_j)$ and suppose $\nu_{Q,l,j}$ exists. Then $\tg(\nu_{Q,l,j}) \notin U_{Q,l,j}$ 
\end{lemma}

\begin{proof}
The proof is similar to the proof of the previous lemma, but simpler. Again, write 
$$ (\om_0, \ldots , \xi_{Q,l,j}, 0, \ldots ,0) := \iota(\nu_{Q,l,j}),$$
such that $\ph_\eta(\sum_i \om_i y^i) = \nu_{Q,l,j}$.
By Definition~\ref{defQlj}, we have $1 < l \leq \ord_Q(E_j+p^nD_j)=p^n$. This gives $\xi_{Q,l,j}$ a term of order $-l$, which is recorded by $g'_{j+1,Q}$ (see Definition~\ref{defg'sec}). To see that this contribution is not in $U_{Q,l,j}$, note that $\ord_Q(\xi_{P,a,i}) > \ord_Q(\xi_{Q,l,j})$, such that the coefficient of $t^{-l}dt_Q$ is zero in $\xi_{P,a,i}$. This finishes the proof that $\tg(\nu_{Q,l,j}) \notin U_{Q,l,j}$.
\end{proof}

\begin{lemma} \label{lemQljQ'}
The elements 
$$ \left\{ \tg_0(\nu_{Q',l,j}) \; : \; (Q',l,j) \in T \text{ and $\nu_{Q',l,j}$ exists} \right\} $$ are linearly independent.
\end{lemma}

\begin{proof}
Again, write
$$ (\om_0, \ldots , \xi_{Q,l,j}, 0, \ldots ,0) := \iota(\nu_{Q,l,j}).$$
By Definition~\ref{defQlj}, we have $1<l\leq p^{n-1}d_{Q'}(p-1-j)$ with $l \not \equiv 1 \bmod p^n$. Recall also from Definition~\ref{defg'global} that $g'_0$ is the natural projection from $\G_0$ to $\G_0/\G_{-1} \cong M_{-1}$, which is only supported at $Q'$. The map $g'_0$ records the coefficients between $t_{Q'}^{p^{n-1}d_{Q'}(p-1-j)}dt_{Q'}$ and $t_{Q'}^{-2}dt_{Q'}$ excluding the exponents that are $-1 \bmod p^n$. We prove linear independence by contradiction: assume that a linear combination of the $\tg_0(\nu_{Q',l,j})$ equals zero in $M_{-1}$, meaning it equals a differential that is regular at $Q'$:
\begin{equation} \label{eqQ'ljlinindep}
    \sum_{(Q',l,j) \in T} c_{Q',l,j} \om_{Q',l,j} = \om' \in H^0(Y, \ker V_Y^n).
\end{equation} 
As in~\cite[Lemma 6.23]{BoCaASc}, we prove that each $c_{Q',l,j}$ is zero by descending induction on $j$. To start with, there are no triples of the form $(Q',l,p-1)$. Now assume that $c_{Q',l,i}=0$ for each $l$ and for $i>j$. Assuming this, we prove that $c_{Q',l,j}=0$ is zero for each $l$. Using $(\tau-1)^j$, we reduce the power of $y$ by $j$ to distill the triples of the form $(Q',l,j)$. Applying $(\tau-1)^j$ to \eqref{eqQ'ljlinindep} gives
$$ (\tau-1)^j \sum_{(Q',l,j) \in T} c_{Q',l,j} \om_{Q',l,j} = j! \sum_{\{ l: (Q',l,j) \in T\}} c_{Q',l,j} \xi_{Q',l,j} = (\tau-1)^j \om'.$$
The terms on the left-hand side all have a different negative valuation at $Q'$, whence the sum can only be regular at $Q'$ if every $c_{Q',l,j}$ is zero. On the other hand, $(\tau-1)^j \om'$ is regular at $Q'$ since $\om'$ is regular at $Q'$. This implies that each $c_{Q',l,j}$ is zero; the elements $\om_{Q',l,j}$ are linearly independent, which finishes the proof.
\end{proof}

Using the preceding three lemmas, we can prove Proposition~\ref{propNTB}.
\begin{proof}[Proof of Proposition~\ref{propNTB}]
We bound both $N_1(X,\pi)$ and $N_2(X,\pi)$ from below. Lemma~\ref{lemQljbrp} and Lemma~\ref{lemQljDj} together imply 
\begin{align*}
    N_1(X,\pi) &= \dim_k \im(\tg) \geq \# \{(Q,l,j) \; : \; \nu_{Q,l,j} \text{ exists and } Q \neq Q' \}.
\intertext{On the other hand, Lemma~\ref{lemQljQ'} implies} 
    N_2(X,\pi) &= \dim_k \im(\tg_0) \geq \# \{(Q',l,j) \; : \; \nu_{Q',l,j} \text{ exists} \}.
    \intertext{Combining these two lower bounds yields}
    N(X,\pi) &= N_1(X,\pi) + N_2(X,\pi) \geq \# \{(Q,l,j) \; : \; \nu_{Q',l,j} \text{ exists}\} = \# T - B,
\end{align*}
as desired.

\end{proof}

\subsection{Explicit bounds}

We have now amassed all the tools needed to prove the main theorem. Denote by $\sigma_p(d,i,n)$ the number of integers $0 < l \leq \lceil id/p \rceil$ that are not $1 \bmod p^n$ and such that the first $n$ digits in the $p$-adic expansion of $(1-l)/d$ sum up to at most $p-1-i$.

\begin{theorem} \label{thmmain}
We have
\begin{equation} \label{equpperbound}
    a_Y^n \leq pa_X^n + \sum_{Q\in S} \sum_{i=1}^{p-1} \left( \left \lfloor \frac{id_Q}{p} \right\rfloor - \left\lfloor \frac{id_Q}{p^{n+1}} \right\rfloor -\sigma_p(d_Q,i,n) \right).
\end{equation}
Moreover, for every $0\leq j \leq p-1$, we have
\begin{equation} \label{eqlowerbound}
    a_Y^n \geq \sum_{Q \in S} \sum_{i=j}^{p-1} \left( \left\lfloor \frac{id_Q}{p} \right\rfloor - \left\lfloor \frac{id_Q}{p} - \left(1-\frac{1}{p^n}\right) \frac{jd_Q}{p} \right\rfloor \right).
\end{equation}
\end{theorem}

\begin{proof}
As before, we have an Artin-Schreier cover $\pi:Y \to X$ and we denote by $S$ the set of points on $X$ where $\pi$ is ramified. For $0 \leq i \leq p-1$, let $E_i$ be the divisors in \eqref{exseqgal} used to split the filtration of sheaves induced by the Galois action. Let $m(X)$ be the Tango number of $X$ and let the divisors $D_i$ be as in Lemma~\ref{lemsecs}. Note that, by adding points if necessary, we can choose the $D_i$ such that $\deg(D_i) > m(X)$, allowing us to apply Tango's theorem and hence Proposition~\ref{proptango}. We write
\begin{align*}
    A_i &:= \sum_{Q\in S'} \left\lfloor \frac{n_{Q,i}}{p^n} \right\rfloor [Q] \\
    R_i &:= E_i-pA_i =: \sum_{Q\in S} r(Q,i) [Q] \\
    s(Q,i)&:= n_{Q,i} - \left\lceil \frac{n_{Q,i}}{p^n} \right\rceil - (p^n-1) \left\lfloor \frac{n_{Q,i}}{p^n} \right\rfloor .
\end{align*}
Note that $s(Q,i)$ equals $r(Q,i)-1$ if $r(Q,i)$ is non-zero and equals zero if $r(Q,i)$ is zero.

We prove the upper bound first. Proposition~\ref{propU(X,pi)} gives
\begin{align} \label{equpbndfirst}
    a_Y^n &= U^n(X,\pi) \nonumber \\
    &:= \sum_{i=0}^{p-1} \dim_k H^0(X,\ker V_X^n(F_*^n(E_i+p^nD_i))) - N(X,\pi) 
\end{align}
We use Corollary~\ref{corwrong} to compute the first term and Proposition~\ref{propNTB} to bound the second term from below. For each $i$, Corollary~\ref{corwrong} gives
\begin{align} \label{equpbndtango}
    \dim_k H^0(X,\ker V_X^n(F_*^n(E_i+p^nD_i)))&= \d(H^0(X,\Om_X^1(p^n(D_i+A_i)+R_i)))\\ \nonumber
    &=(p^n-1)\deg(D_i+A_i) + \sum_{Q\in S} s(Q,i) \\ \nonumber
    = \frac{(p^n-1)d_{Q'}(p-1-i)}{p} + \sum_{Q\in \supp(D_i)} &(p^n-1) + \sum_{Q \in S} \left((p^n-1) \left\lfloor \frac{n_{Q,i}}{p^n} \right\rfloor +s(Q,i)  \right).
\end{align} 
Since the divisor $D_i+A_i$ is supported on $S_i'$, we have distributed $\deg(A_i+D_i)$ over the partition $S_i=\{Q'\} \cup \supp(D_i) \cup S$. Now, Proposition~\ref{propNTB} states $N(X,\pi) \geq \#T-B$. We first note that by the proof of Corollary~\ref{corwrong} we have $B \leq pa_X^n$. We now count the number of triples $(Q,l,i)$ in $T$. For $Q=Q'$ there are $ \frac{(p^n-1)d_{Q'}(p-1-i)}{p} $ such triples, since $\ord_{Q'}(E_i+p^nD_i)=p^{n-1}d_{Q'}(p-1-i)$. Note that this cancels the contribution from $Q'$ in \eqref{equpbndtango}. For $Q\in \supp(D_i)$, there are $p^n-1$ such triples, which also cancels the contribution in \eqref{equpbndtango}. Finally, for $Q\in S$, by Definition~\ref{defQlj} and Lemma~\ref{lemQljbrp}, the number of triples is the number of integers $0 < l\leq n_{Q,i}$ that are not $1 \bmod p^n$ and such that the first $n$ digits in the $p$-adic expansion of $\frac{1-l}{d_Q}$ sum up to at most $i$. This is by definition $\sigma_p(d_Q,p-1-i,n)$. Substituting all of this into \eqref{equpbndfirst} yields
\begin{align*}
    a_Y^n &\leq pa_X^n + \sum_{i=0}^{p-1} \sum_{Q\in S} \left((p^n-1) \left\lfloor \frac{n_{Q,i}}{p^n} \right\rfloor + s(Q,i) - \sigma_p(d_Q,p-1-i,n) \right) \\
    &= pa_X^n + \sum_{Q\in S} \sum_{i=0}^{p-1} \left( n_{Q,i} - \left\lceil \frac{n_{Q,i}}{p^n} \right\rceil - \sigma_p(d_Q,p-1-i,n)\right) \\
    &= pa_X^n + \sum_{Q\in S} \sum_{i=0}^{p-1} \left( \left\lceil \frac{(p-1-i)d_Q}{p} \right\rceil - \left\lceil \frac{(p-1-i)d_Q}{p^{n+1}} \right\rceil - \sigma_p(d_Q,p-1-i,n) \right).
\end{align*}
The expression in \eqref{equpperbound} is then obtained by using $\lceil x \rceil = \lfloor x \rfloor+1$ for $x \notin \Z$, realizing the summand is zero for $i=p-1$ and then re-indexing $i \mapsto p-1-i$.

We now prove the lower bound. Proposition~\ref{propL(X,pi)} yields that for every $j$ we have
\begin{equation} \label{eqlowbndfirst}
    a_Y^n \geq \dim_k U_j - \sum_{i=0}^j \sum_{Q\in S_i'} c(i,j,Q)
\end{equation}
By Lemma~\ref{lemUj} and Corollary~\ref{corwrong}, we can compute the first term:
$$\dim_k U_j = \sum_{i=0}^j \d (H^0(X,\Om_X^1(E_i+p^nD_i))) $$ 
$$= \sum_{i=0}^j \left( \frac{(p^n-1)d_{Q'}(p-1-i)}{p} + \sum_{Q\in \supp(D_i)} (p^n-1) + \sum_{Q \in S} \left((p^n-1) \left\lfloor \frac{n_{Q,i}}{p^n} \right\rfloor +s(Q,i)  \right) \right).$$
Recall that by Definition~\ref{defL(X,pi)} of $c(i,j,Q)$, the contributions from $Q'$ and from $Q \in \supp(D_i)$ are canceled in \eqref{eqlowbndfirst}. For $Q\in S$, $c(i,j,Q)$ equals the number of integers $m$ that are $-1 \bmod p^n$ such that $-n_{Q,i}-d_Qp^{n-1}(j-i)\leq m \leq -n_{Q,i}$. More explicitly,
$$c(i,j,Q)= \left\lceil \frac{n_{Q,j} - d_Qp^{n-1}(j-i) }{p^n} \right\rceil - \left\lceil \frac{n_{Q,i}}{p^n} \right\rceil $$
What remains is
\begin{align*}
    a_Y^n &\geq \sum_{i=0}^j \sum_{Q\in S} \left( (p^n-1) \left\lfloor \frac{n_{Q,i}}{p^n} \right\rfloor +s(Q,i) - \left\lceil \frac{n_{Q,j} - d_Qp^{n-1}(j-i) }{p^n} \right\rceil + \left\lceil \frac{n_{Q,i}}{p^n} \right\rceil \right) \\
    &= \sum_{Q\in S} \sum_{i=0}^j \left( n_{Q,i} -  \left\lceil \frac{\frac{(p-1-j)d_Q}{p} - d_Qp^{n-1}(j-i) }{p^n} \right\rceil \right) \\
    &= \sum_{Q\in S} \sum_{i=0}^j \left( \left\lceil \frac{(p-1-i)d_Q}{p} \right\rceil - \left\lceil \frac{(p-1-i)d_Q}{p} - \left(1-\frac{1}{p^n}\right)\frac{(p-l-j)d_Q}{p} \right\rceil \right).
\end{align*}
The expression in \eqref{eqlowerbound} is obtained by using $\lceil x \rceil = \lfloor x \rfloor+1$ for $x \notin \Z$ and re-indexing $i \mapsto p-1-i$.
\end{proof}

As an application of Theorem~\ref{thmmain}, we can generalize~\cite[Lemma 3.3]{Prieshyp2}, replacing the projective line by any ordinary base curve.

\begin{corollary} \label{corp=2}
If $p=2$ and $X$ is ordinary, then the bounds in~\ref{thmmain} coincide: 
$$a_Y^n = \sum_{Q \in S} \frac{d_Q-1}{2} - \left \lfloor \frac{d_Q}{2^{n+1}} \right \rfloor.$$
In particular, the ranks of $V_Y^n$ are determined by the ramification invariants.
\end{corollary}

\begin{proof}
We study the upper bound first. For each $Q\in S$, the sum has only the term corresponding to $i=1$. One notes $\sigma_2(d_Q,i,n)=0$, such that 
$$ a_Y^n \leq \sum_{Q\in S} \left \lfloor \frac{d_Q}{2} \right\rfloor - \left \lfloor \frac{d_Q}{2^{n+1}} \right \rfloor. $$
Since $d_Q$ is odd, $\lfloor d_Q/2\rfloor = (d_Q-1)/2$.

For the lower bound, we see that the expression is zero when $j=0$. Setting $j=1$ gives
\begin{align*}
    a_Y^n &\geq  \sum_{Q\in S} \left\lfloor \frac{d_Q}{2} \right\rfloor- \left\lfloor \frac{d_Q}{2}- \left(1-\frac{1}{2^n} \right) \frac{d_Q}{2} \right\rfloor \\
    &= \sum_{Q\in S} \frac{d_Q-1}{2} - \left \lfloor \frac{d_Q}{2^{n+1}} \right \rfloor.
\end{align*}
If $X$ is ordinary, then $V_X$ is bijective, so that $a_X^n=0$ for every $n$. By the argument above, this implies that every $a_Y^n$ is fixed by the ramification data.
\end{proof}

If $X$ is the projective line, the same formula for $a_Y^n$ can be derived by following the proof of~\cite[Proposition 3.4]{Prieshyp2}.

\begin{remark} \label{rmkaYn>aYn-1}
Forgetting specific properties of the Cartier operator, basic linear algebra tells us $a_Y^n \geq a_Y^{n-1}$. In fact, equality holds if and only if $a_Y^{n-1}=g_Y-s_Y$ (recall that $s_Y$ is the $p$-rank). Recall that $g_Y$ and $s_Y$ can be computed using the Riemann-Hurwitz formula and the Deuring-Shafarevich formula, respectively. Thus if $a_Y^{n-1}<g_Y-s_Y$, we obtain the bound 
$$a_Y^n \geq a_Y^{n-1} +1 \geq L^{n-1}(X,\pi)+1 $$ 
or an iteration thereof. In some cases, especially when the $d_Q$'s are small relative to $p^n$, this bound is better than the lower bound of Theorem~\ref{thmmain}.
\end{remark}

\subsection{Sharpness of the upper bound}

We prove that the upper bound is sharp in the case $X=\PP^1$ and $\#S=1$. In particular, we show that for any $p \nmid d$, the curve 
$$Y: y^p-y= t^d$$
realizes the upper bound for $a_Y^n$ for all $n>0$. In the computation of the upper bound, only the leading term of $f$ is used, so it is not surprising that the upper bound is achieved when $f$ is in fact a monomial.  

We first compute $a_Y^n$, loosely following the strategy in \cite[Lemma 3.1]{Pries05}. Given integers $b$ and $n$, define the integer $h_{b,n}$ as follows. If the $p$-adic expansion of $(-1-b)/d$ is given by
$$(-1-b)/d = c_0 + c_1 p + \ldots + c_{n-1} p^{n-1} + \ldots,$$
set $h_{b,n} := \sum_{i=0}^{n-1} c_i$.

\begin{proposition} \label{propf=t^d}
Let $Y$ be given by $y^p-y=t^d$. Then we have
$$a_Y^n = \sum_{b=0}^{d-2} \min \left \{h_{b,n} \; , \; p-\lceil (p+1+bp)/d\rceil \right\}.$$
\end{proposition}
\begin{proof}
By \cite[Lemma 1]{Sulli}, a basis of $H^0(Y, \Om_Y^1)$ is given by differentials $y^r t^b dt$ satisfying $r,b \geq 0$ and
\begin{equation*}
rd+bp \leq (p-1)(d-1)-2.
\end{equation*}
We apply $V_Y^n$ to such a basis differential using Lemma~\ref{lemVYinVX}. This method is comparable to the proof of Lemma~\ref{lemQljbrp}. By substituting $f=t^d$, the operator $W_{l,m}$ simplifies to
$$W_{l,m} (\omega) = V_X \left( \binom{l}{m} (-f)^{l-m} \omega \right) = (-1)^{l-m} \binom{l}{m} V_X \left( t^{d(l-m)} \omega \right).$$
Lemma~\ref{lemVYinVX} yields
$$V_Y^n(y^r t^b dt) = \sum_{j=0}^{p-1} \left(
\sum_{j\leq j_{1}\leq \ldots \leq r}
W_{j_{1},j}( W_{j_{2},j_{1}}( \ldots W_{r,j_{n-1}}(t^b dt) \ldots )) \right) y^j.$$
We begin in the innermost parenthesis by computing 
$$W_{r,j_{n-1}}(t^b dt) = (-1)^{r-j_{n-1}} \binom{r}{j_{n-1}} V_X \left(t^{d(r-j_{n-1})+b} dt \right).$$
This is non-zero if and only if $d(r-j_{n-1})+b \equiv -1 \bmod p$. This in turn is equivalent to $r-j_{n-1} \equiv (-1-b)d^{-1} \bmod p$. Since we have $0 \leq j_{n-1} \leq p-1$, the solution to the congruence is unique. Applying $W_{j_{n-1},j_{n-2}}$ and repeating this inductively shows that $V_Y^n(y^r t^b dt)$ has a non-zero term if and only if
$$p^n \mid b+1+d_Q\left(r-j_{n-1} + (j_{n-1} - j_{n-2})p + \ldots + (j_1-j)p^{n-1} \right).$$
It then follows that $j=r-h_{b,n}$, so that $j \geq 0$ is equivalent to $r \geq h_{b,n}$. Thus we obtain
\begin{equation*}
V_Y^n(y^r t^b dt) \begin{cases}
= 0 & \hbox{if $r<h_{b,n}$} \\
\in \spn_k \{y^{r-h_{b,n}} t^{(h_{b,n}d+b+1)/p^n-1} dt\} &\hbox{if $r \geq h_{b,n}$}.
\end{cases}
\end{equation*}
Since the differentials $y^{r-h_{b,n}} t^{(h_{b,n}d+b+1)/p^n-1} dt$ are distinct for distinct pairs $(b,r)$, it follows that $a_Y^n$ is given by
\begin{align*}
a_Y^n &= \# \left\{ (b,r) \; | \; r,b \geq 0 \; , \; rd+bp \leq (p-1)(d-1)-2, \; r<h_{b,n}  \right\} \\
&= \sum_{b=0}^{d-2} \min \left \{h_{b,n} \; , \; p-\lceil (p+1+bp)/d\rceil \right\},
\end{align*}
which is what we wanted to prove.
\end{proof}

We now prove that this expression equals the upper bound from Theorem~\ref{thmmain}.

\begin{corollary} \label{corsharp}
Let $Y$ be an Artin-Schreier curve given by $y^p-y=t^d$ and let $n$ be a natural number. Then $Y$ attains the upper bound of Theorem~\ref{thmmain}:
$$a_Y^n = \sum_{i=1}^{p-1} \left( \left \lfloor \frac{id}{p} \right\rfloor - \left\lfloor \frac{id}{p^{n+1}} \right\rfloor -\sigma_p(d,i,n) \right).$$
\end{corollary}
\begin{proof}
We've seen in Proposition~\ref{propf=t^d} that $a_Y^n$ equals the cardinality of the set 
$$A := \left\{ (b,r) \; | \; r,b \geq 0 \; , \; rd+bp \leq (p-1)(d-1)-2, \; r<h_{b,n}  \right\}.$$
On the other hand, the proof of Theorem~\ref{thmmain} gives the following expression for the upper bound of $a_Y^n$:
$$\sum_{i=0}^{p-1} \left( n_{Q,i} - \lceil n_{Q,i}/p^n \rceil - \sigma_p(d,p-1-i,n) \right).$$
This equals the cardinality of the set $B$ consisting of pairs $(i,l)$ with the following properties:
\begin{enumerate}
    \item $0 \leq i \leq p-1$; 
    \item $2 \leq l \leq n_{Q,i}$;
    \item $l \not \equiv 1 \bmod p^n$;
    \item $i$ is smaller than the sum of the first $n$ digits in the $p$-adic expansion of $(1-l)/d$.
\end{enumerate}
Using the definition $n_{Q,i}=\lceil (p-1-i)/p \rceil$, one verifies that the first two properties together are equivalent to the equation
$$id + (l-2)p \leq (p-1)(d-1)-2.$$
Furthermore, the fourth property is equivalent to $i<h_{l-2,n}$ and then the third property is automatically satisfied. Thus the sets $A$ and $B$ are in bijection via $(b,r)=(l-2,i)$. We conclude that indeed $a_Y^n$ coincides with the upper bound of Theorem~\ref{thmmain}.
\end{proof}

\section{Examples} \label{sec7}

In this section we apply the bounds for $a_Y^n$ proved in Theorem~\ref{thmmain} to a few examples. In doing so we study how the bounds behave as a function of $n$. The space $H^0(Y, \Om_Y^1)$ decomposes as
$$ H^0(Y, \Om_Y^1) = H^0(Y, \Om_Y^1)^{\textup{bij}} \oplus H^0(Y, \Om_Y^1)^{\textup{nil}},$$
where $V_Y$ acts bijectively on $H^0(Y, \Om_Y^1)^{\textup{bij}}$ and nilpotently on $H^0(Y, \Om_Y^1)^{\textup{nil}}$. The dimension of $H^0(Y, \Om_Y^1)^{\textup{bij}}$ equals the $p$-rank $s_Y$, so the ranks of higher powers of $V_Y^n$ approach $s_Y$. Therefore one would expect that the bounds get closer and closer to $g_Y-s_Y=\dim \ker(V_Y^{g_Y})$ as $n$ increases. Note that combining the Riemann-Hurwitz formula and the Deuring-Shafarevich formula yields
\begin{equation*} \label{eqg-s}
g_Y-s_Y = p(g_X-s_X) + \sum_{Q \in S} \frac{(p-1)(d_Q-1)}{2}.
\end{equation*} 
By Remark~\ref{rmkaYn>aYn-1}, the lower bound reaches $g_Y-s_Y$ for some $n\leq g_Y-s_Y$, but generally it is reached faster than Remark~\ref{rmkaYn>aYn-1} prescribes. Looking at the upper bound $U^n(X,\pi)$, we see that $\sigma_p(d,i,n)$ vanishes for sufficiently large $n$. Then the upper bound converges to 
$$ pa_X^n+ \sum_{Q \in S} \sum_{i=1}^{p-1} \left \lfloor \frac{id_Q}{p} \right \rfloor = p(g_X-s_X)+ \sum_{Q \in S} \frac{(p-1)(d_Q-1)}{2}= g_Y-s_Y.$$
In line with this expectation, one would expect both bounds to be increasing as a function of $n$, with the distance between the bounds eventually decreasing.

Computations that could not be done by hand have been done using Magma (\cite{MAGMA}).

\begin{example} \label{egp=3d=100}
This example demonstrates how the bounds can behave as a function of $n$. Let $p=3$ and let $Y \to \PP^1$ be branched at a single point $Q$ with ramification invariant $d_Q=100$. This curve has genus $99$ and $p$-rank $0$. The bounds of Theorem~\ref{thmmain} and Remark~\ref{rmkaYn>aYn-1} are displayed in Table~\ref{tabp=3d=100} and in Figure~\ref{plotp=3d=100}.

\begin{table}[ht] 
\begin{tabular}{l|llllllllll}
$n$          & $1$  & $2$  & $3$  & $4$  & $5$  & $6$  & $7$  & $8$  & $9$  & $10$ \\ \hline
$L^n(X,\pi)$ & $44$ & $59$ & $64$ & $66$ & $67$ & $68$ & $69$ & $70$ & $71$ & $72$ \\
$U^n(X,\pi)$ & $55$ & $82$ & $93$ & $98$ & $99$ & $99$ & $99$ & $99$ & $99$ & $99$
\end{tabular}
\caption{The bounds $L^n(X,\pi)$ and $U^n(X,\pi)$ for $X=\PP^1$, $S=\{Q\}$, $d=100$ and varying $n$.}
\label{tabp=3d=100}
\end{table}

\begin{figure} 
\begin{tikzpicture} 
\begin{axis}[
    xlabel={$n$},
    ylabel={$a_Y^n$},
    xmin=1, xmax=10,
    ymin=0, ymax=100,
    xtick={1,2,3,4,5,6,7,8,9,10},
    ytick={0,20,40,60,80,100},
    legend pos=south east,
    ymajorgrids=true,
]

\addplot+[
    only marks,
    color=blue,
    mark=*,
    ]
    coordinates {
    (1,44)(2,59)(3,64)(4,66)(5,67)(6,68)(7,69)(8,70)(9,71)(10,72)
    };
    
    \addplot+[
    only marks,
    color=red,
    mark=*,
    ]
    coordinates {
    (1,55)(2,82)(3,93)(4,98)(5,99)(6,99)(7,99)(8,99)(9,99)(10,99)
    };
\end{axis}
\end{tikzpicture}
\caption{A plot of values of the bounds $L^n(X,\pi)$ and $U^n(X,\pi)$ for Artin-Schreier curves in characteristic $3$ branched at a single point with ramification invariant $100$.}
\label{plotp=3d=100}
\end{figure}
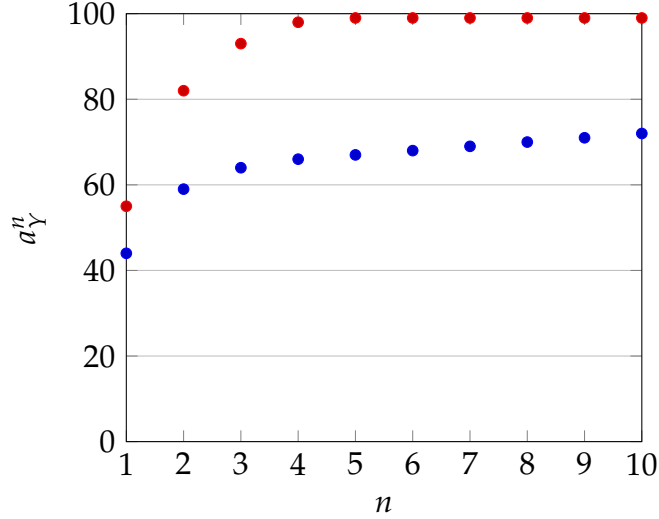

The lower bound of Theorem~\ref{thmmain} stabilizes at $66$, so for $n\geq 5$ we apply Remark~\ref{rmkaYn>aYn-1}. For $n\geq 5$, the upper bound is equal to the genus $g_Y$, as the curve has $p$-rank zero. As a consequence of Corollary~\ref{corsharp}, the Artin-Schreier curve
$$Y: y^3-y=t^{100} $$ simultaneously realizes all the upper bounds: $a_Y^n=U^n(X,\pi)$ for all $n\geq 1$. Whether the lower bounds are realized (simultaneously) by a curve is a more difficult question and a potential research topic for the future, perhaps following an approach like the one in~\cite{11a}.

When we set $d=101$, the behavior is very similar, suggesting that the residue class of $d$ modulo $p$ does not exert a large influence.
\end{example}

\begin{example}
We generalize Example~\ref{eg7}. Let $Y$ be an Artin-Schreier curve given by $$Y: y^7-y=t^{-4}+c_3t^{-3}+c_2t^{-2} + c_1 t^{-1}.$$ Theorem~\ref{thmmain} prescribes the bounds $6 \leq a_Y^2 \leq 9$. However, setting $n=1$ in~\ref{thmmain} gives $6 \geq a_Y^1$, such that Remark~\ref{rmkaYn>aYn-1} yields $a_Y^2 \geq 7$.

As $X$ is the projective line, we do not need a point $Q'$. This implies an equality that $H^0(Y,\ker V_Y^2) = H^0(X,\G_0)$. In Example~\ref{eg7} we saw that $\bigoplus_{i=0}^{p-1} H^0(\PP^1,\ker V_{\PP^1}^2(F_*^2E_i))$ is spanned by the elements $\nu_{j,l}$, with $2\leq l \leq n_{Q,j},$ that have $t^{-l}$ in the $j$-th position and zero elsewhere. It is clear from Definition~\ref{defGj} and the definition of $\ph^{-1}$ in Proposition~\ref{propphieta} that $\nu_{0,2}$, $\nu_{0,3}$ and $\nu_{0,4}$ lie in $H^0(X,\G_0)$. Similarly we conclude that $\nu_{1,2}, \nu_{1,3} \in H^0(X,\G_0)$, since $$V_{\PP^1}(t^{-l}dt)=V_{\PP^1}((-f)t^{-l}dt)=0$$
for $l \in \{2,3\}$. As moreover $V_{\PP^1}^2((-f)^2t^{-l}dt)=0$ for $l \in \{2,3\}$, we deduce that $\nu_{2,2}$ and  $\nu_{2,3}$ lie in $H^0(X,\G_0)$. For $\nu_{3,2}$, we compute
\begin{align*}
    V_{\PP^1}^2((-f)^3t^{-2}dt)&=0 \\
    V_{\PP^1}((-f)^2t^{-2}dt)&=(c_3^2+2c_2)t^{-2}dt \\
    V_{\PP^1}((-f) V_{\PP^1}((-f)^2 t^{-2}dt))&=0.
\end{align*}
This implies $\nu_{3,2} \in H^0(X,\G_0)$.

Now, if we compute $\om_0$ for $\nu_{4,2}$, following Proposition~\ref{propphieta}, then the only possibly non-zero term is 
$$s_0\left(V_{\PP^1}((-f)^2 V_{\PP^1}((-f)^2 t^{-2} dt  \right)=(c_3^2+2c_2)^2 t^{-50} dt.$$
Thus $a_Y^2=9$ if $c_3^2+2c_2=0$ and $a_Y^2=8$ otherwise. This is an example where the lower bound for $a_Y^2$ is not sharp.
\end{example}

\section*{Data availabilty statement}

Some minor programming in Magma (\cite{MAGMA}) was involved in Example~\ref{egp=3d=100}. The code is available at \url{https://warwick.ac.uk/fac/sci/maths/people/staff/groen/powers_cartier_code/}.

\section*{Conflict of interest statement}

The author is funded by the University of Warwick and hereby declares that he has no conflicts of interest to disclose.

\section*{University of Warwick Rights Retention Statement}

For the purpose of open access, the author has applied a Creative Commons Attribution (CC-BY) licence to any Author Accepted Manuscript version arising from this submission.

\bibliographystyle{plain}
\bibliography{main}

\end{document}